\documentclass[twoside, onecolumn, 10pt]{article}
\usepackage{latexsym}
\usepackage{amssymb}
\usepackage{mathrsfs}
\usepackage{amsmath}
\usepackage{marvosym}
\newtheorem{theorem}{Theorem}[section]
\newtheorem{lemma}[theorem]{Lemma}
\newtheorem{corollary}[theorem]{Corollary}

\newtheorem{problem}[theorem]{Problem}
\newtheorem{conjecture}[theorem]{Conjecture}

\newcommand{\nb}{\nonumber}
\newcommand{\R}{{\mathscr R}}

\def\bc{\begin{center}}
\def\ec{\end{center}}

\def\nb{\nonumber}

\begin{document}

\begin{center}
{\LARGE Analytical solutions to some optimization problems
on ranks and inertias of matrix-valued functions subject to linear matrix inequalities}
\end{center}

\begin{center}
{\Large  Yongge Tian}
\end{center}

\begin{center}
{\em {\footnotesize CEMA, Central University of Finance and Economics,
 Beijing 100081, China}}
\end{center}

\renewcommand{\thefootnote}{\fnsymbol{footnote}}
\footnotetext{ {\it  {\rm E-mail Address}: yongge.tian@gmail.com}}

\noindent {\bf Abstract.} \ {\small
Matrix rank and inertia optimization problems are a class of discontinuous optimization problems,
in which the decision variables are matrices running over certain feasible matrix sets, while
the ranks and inertias of the variable matrices are taken as integer-valued objective functions.
In this paper, we establish a group of explicit formulas for calculating the maximal and minimal values
of the rank- and inertia-objective functions of the Hermitian matrix  expression $A_1 - B_1XB_1^{*}$
 subject to the linear matrix inequality  $B_2XB_2^{*} \succcurlyeq  A_2$ $(B_2XB_2^{*} \preccurlyeq  A_2)$
 in the L\"owner partial ordering, and give applications of these formulas in characterizing behaviors of
 some constrained matrix-valued functions.\\

\noindent {\em Key Words}:  Matrix-valued function;  matrix equation; LMI;
rank; inertia;  integer-valued objective function; feasible matrix set; generalized inverses of matrices;
 optimization; L\"owner partial ordering

\medskip

\noindent {\em Mathematics Subject Classifications:} 15A24; 15A39; 15A45; 15B57; 49K30; 65K10; 90C11; 90C22}

\section{Introduction}

\renewcommand{\theequation}{\thesection.\arabic{equation}}
\setcounter{section}{1}
\setcounter{equation}{0}

Throughout this paper,
\begin{enumerate}
\item[]  ${\mathbb C}^{m\times n}$ stands for the set of all $m\times n$ complex
matrices;\

\item[] ${\mathbb C}_{{\rm H}}^{m}$
stands for the set of all $m\times m$ complex Hermitian matrices;

\item[] $A^{*}$,  $r(A)$ and ${\mathscr R}(A)$ stand
for the conjugate transpose, rank and range (column space) of a
matrix $A\in {\mathbb C}^{m\times n}$, respectively;

\item[]  $I_m$ denotes the identity matrix of order $m$;

\item[]  $[\, A, \, B\,]$ denotes a row
block matrix consisting of $A$ and $B$;

\item[]  the Moore--Penrose inverse of $A\in {\mathbb
C}^{m\times n}$, denoted by $A^{\dag}$, is defined to be the unique
solution $X$ satisfying the four matrix equations $AXA = A,$ $XAX =
X,$ $(AX)^{*} = AX$ and $(XA)^{*} =XA$;

\item[] the symbols $E_A$ and $F_A$ stand for $E_A = I_m - AA^{\dag}$ and $F_A = I_n-A^{\dag}A$;

\item[] $i_{+}(A)$ and $i_{-}(A)$, called the partial inertia of  $A \in {\mathbb C}_{{\rm H}}^{m}$,
 are defined to be  the numbers of the positive and negative eigenvalues of $A$ counted with
multiplicities, respectively;

\item[]  $A \succcurlyeq  0$ ($A \succ  0$) means that $A$ is Hermitian positive semi-definite (Hermitian positive definite);

\item[]  two $A,  \, B \in {\mathbb C}_{{\rm H}}^{m}$ are said to satisfy the inequality $A \succcurlyeq  B$
($A \succ  B)$ in the L\"owner partial ordering if $A - B$ is positive semi-definite (positive definite).
\end{enumerate}

The matrix approximation problem is to approximate optimally, with respect to some criteria,
a matrix by one of the same dimension from a given feasible matrix set. Assume that $A$ is a matrix to
be approximated. Then a conventional statement of general matrix optimization  problems of $A$
from this point of view can be written as
\begin{align}
{\rm minimize}\ \rho(\, A - Z \,) \ \  {\rm subject \ to}  \ Z  \in {\cal S},
\label{ww11}
\end{align}
where $\rho(\cdot)$ is certain objective function, which is usually taken as the determinant, trace,
 norms, rank, inertia of  matrix, and ${\cal S}$  is a given feasible matrix set.
 A best-known case of (\ref{ww11})  is to  minimize the norm $ \|\, A - Z \,\|^2_{F}$ subject to  $Z  \in {\cal S}$.

In this paper, we assume that the objective function $\rho(\cdot)$ in (\ref{ww11}) is taken as the rank or inertia of matrix,
and $A \in {\mathbb C}^{m \times m}$ is a Hermitian matrix. The rank and inertia of matrix, as  objective functions, are often used when finding feasible matrices $Z$ such that
resulting $A - Z$ attains its maximal possible rank or inertia (is nonsingular or definite when square), or finding feasible
matrix $Z$ such that $A - Z$ attains the minimal rank or inertia as possible (called low-rank or low-inertia matrix completion).
This kind of  problems are usually called the matrix rank-optimization and inertia-optimization problems, or
 matrix rank and inertia completion problems in the literature. Generally speaking, matrix rank and inertia optimization
 problems are a class of discontinuous optimization problems, in which the decision variables are matrices running over
 certain matrix sets, while the ranks and  inertias of the variable matrices are taken as integer-valued objective
 functions. This kind of optimization problems can generally be written as
\begin{align}
& \text{maximize}\ r(\, A - Z\,) \ \ \ \ \ \,   {\rm subject \ to}  \ \ \  Z \in {\cal S},
\label{ww13}
\\
& \text{minimize}\ r(\, A - Z\,) \ \ \ \ \ \   {\rm subject \ to}  \ \ \  Z \in {\cal S},
\label{ww14}
\\
& \text{maximize} \ i_{\pm}(\, A - Z\,) \ \ \  {\rm subject \ to}  \ \ \  Z \in {\cal S},
\label{ww15}
\\
& \text{minimize} \ i_{\pm}(\, A - Z\,) \ \ \ \ {\rm subject \ to}  \ \ \  Z \in {\cal S},
\label{ww16}
\end{align}
respectively.

The rank and inertia of a Hermitian matrix are two generic concepts in matrix theory for describing the dimension
of the row or column vector space and the sign distribution of the eigenvalues of the matrix, which are well
understood and are easy to compute by the well-known elementary or congruent matrix operations. These two quantities
play an essential role in characterizing relations between two matrices and algebraic properties of matrices.
These two integer indices occur only in finite-dimensional algebras and are not replaceable and cannot be approximated
by other continuous quantities. Because the rank and inertia of a matrix are always
finite nonnegative integers less than or equal to the dimensions of the matrix, it is not hard to give upper
and lower bounds for ranks and inertias of matrices, and the global
maximal and minimal values of the integer-valued objective functions always
exist, no matter what the decision domain ${\cal S}$ is given. Also, due to the integer property of
rank and inertia, inexact or approximate values of maximal and minimal ranks and inertias are less valuable,
so that no approximation methods are allowed to use when finding the maximal and minimal possible ranks and inertias
 of a matrix-valued function. This fact means that solving methods of matrix rank and inertia optimization problems
 are not consistent with any of the ordinary continuous and discrete problems in optimization theory.
It has been known that matrix rank optimization problems are NP-hard in general due to the discontinuity and
combinational nature of rank of a matrix and the algebraic structure of ${\cal S}$.   However, it is really
lucky that we can establish analytical formulas for calculating the extremal ranks of matrix-valued functions
for some special feasible matrix sets ${\cal S}$ by using various expansion formulas for ranks and inertias of matrices
and some tricky matrix operations.

Because  the rank of a matrix can only take finite integers between 0 and the dimensions of the matrix, it is really
expected to establish certain analytical formulas for calculating the maximal and minimal ranks for curiosity.
In recent years, maximization and minimization problems on ranks and inertias of  matrices attract much
attention from both theoretical and practical points of view.  In this paper, we assume that
 $A_i \in \mathbb C_{{\rm H}}^{m_i}$ and $B_i \in \mathbb C^{m_i \times n}$, $ i =1,\, 2$ are given matrices,
and the feasible matrix set ${\cal S}$ is
\begin{align}
{\cal S}_1 = \{ Z = B_1XC_1  \ | X \in  {\mathbb C}_{{\rm H}}^{n}, \ B_2XB^{*}_2 \succcurlyeq  A_2  \} \
{\rm or} \  {\cal S}_2 = \{ Z= B_1XC_1  \ | X \in  {\mathbb C}_{{\rm H}}^{n}, \ B_2XB^{*}_2 \preccurlyeq  A_2  \}.
\label{ww17}
\end{align}
Then, the difference $A_1 - Z$ can be written as the following linear matrix-valued function
\begin{align}
\phi(X) = A_1  - B_1XC_1.
\label{ww18}
\end{align}
The LMIs in (\ref{ww17}), the simplest cases of all  LMIs, could be regarded as extensions of the usual inequalities $bx \geqslant a$ and $bx \leqslant a$ for real numbers.

Under such a formulation, this paper aims at solving the following inequality-constrained matrix optimization problems:

\begin{problem} \label{TS11}
{\rm
For the function  in (\ref{ww18}) and the feasible matrix sets in (\ref{ww17}),
establish explicit formulas for calculating the following extremal ranks and inertias
\begin{align}
 & \max r(\, A_1-B_1XB_1^{*} \,)  \ \ \ \ \ \,  {\rm s.t.}  \ \ X\in {\cal S}_i,  \ \ i = 1,\, 2,
\label{ww19}
\\
& \min r(\, A_1-B_1XB_1^{*} \,)  \ \  \ \ \ \ {\rm s.t.}  \ \ X\in {\cal S}_i, \ \  i = 1,\, 2,
\label{ww110}
\\
 &  \max i_{\pm}(\, A_1-B_1XB_1^{*} \,)  \ \ \ \  {\rm s.t.}  \ \ X\in {\cal S}_i, \ \ i = 1,\, 2,
\label{ww111}
\\
& \min\!i_{\pm}(\, A_1-B_1XB_1^{*} \,) \ \  \ \ \ {\rm s.t.}  \ \ X\in {\cal S}_i, \ \ i = 1,\, 2,
 \label{ww112}
\end{align}
respectively.
}
\end{problem}

\begin{problem} \label{TS12}
{\rm Establish necessary sufficient conditions for the following two linear matrix inequalities (LMIs)
\begin{eqnarray}
 B_1XB^{*}_1 \succcurlyeq  A_1  \ {\rm and} \ B_2XB^{*}_2 \succcurlyeq   A_2  \ \
 (\, B_1XB^{*}_1 \preccurlyeq  A_1  \ {\rm and} \ B_2XB^{*}_2 \preccurlyeq   A_2 \,)
\label{qq15}
\end{eqnarray}
to have a common Hermitian solution and give their common solutions$.$
}
\end{problem}

\begin{problem} \label{TS13}
{\rm For $\phi(X)$ in (\ref{ww18}), establish necessary and
sufficient conditions for the existence of  $\widehat{X}, \, \widetilde{X} \in {\mathbb C}_{{\rm H}}^{n}$
such that
\begin{align}
\phi(\widehat{X}) \preccurlyeq  \phi(X) \preccurlyeq  \phi(\widetilde{X}) \ \ {\rm for \ all} \ \   B_2XB^{*}_2  \succcurlyeq
A_2  \ {\rm and} \  X \in  {\mathbb C}_{{\rm H}}^{n},
\label{qq16}
\\
\phi(\widehat{X}) \preccurlyeq  \phi(X) \preccurlyeq  \phi(\widetilde{X}) \ \ {\rm for \ all} \ \   B_2XB^{*}_2  \preccurlyeq
A_2  \ {\rm and} \  X \in  {\mathbb C}_{{\rm H}}^{n}
\label{qq16a}
\end{align}
hold, respectively, and find analytical expressions of $\widehat{X}$ and $\widetilde{X}$.
}
\end{problem}

The matrix function $\phi(X) = A - BXB^{*}$, as one of the simplest cases among
all matrix maps  with symmetric patterns, attracted much attention in the recent decade,
and many problems on  $\phi(X)$, were considered in the literature.  Some recent work on the
matrix function is summarized below:
\begin{enumerate}
\item[{\rm (i)}] Expansion formulas for calculating the (global extremal) rank and inertia of
$\phi(X)$ when $X$ running over ${\mathbb C}_{{\rm H}}^{n}$,   \cite{LT-nla12,T-laa10,TL}.

\item[{\rm (ii)}] Nonsingularity, positive definiteness, rank and inertia invariance, etc., of
$\phi(X)$, \cite{T-laa10,TL}.

\item[{\rm (iii)}]  Canonical forms of  $\phi(X)$ under generalized singular value decompositions and their algebraic properties, \cite{LT-nla12}.

\item[{\rm (iv)}]  Solutions and least-squares solutions of the matrix equation $\phi(X) = 0$ and their algebraic properties, \cite{LT-jamc,LT-nla12,LTT,T-laa11,T-med,T-nla}.

\item[{\rm (v)}] Minimization of ${\rm tr}[\,\phi(X) \phi^{*}(X)\,]$
s.t. $r[\phi(X)]= \min$,  \cite{T-nla}. \\[-4mm]

\item[{\rm (vi)}]  Solutions of the  matrix inequalities  $\phi(X) \succ  \, (\succcurlyeq, \, \prec , \, \preccurlyeq )\, 0$ and their properties,  \cite{T-laa10,T-mcm}.

\item[{\rm (vii)}]  Formulas for calculating the extremal rank and inertia of
 $\phi(X)$ under the restrictions $r(X) \leqslant k$ and/or $\pm
X \succcurlyeq 0$, \cite{T-laa10,T-na,T-mcm,TL}.

\item[{\rm (viii)}]  Formulas for calculating the extremal rank and inertia of
 $\phi(X)$ subject to the Hermitian solution of a consistent
matrix equation $CXC^{*} = D$, \cite{LT-jota}.

\item[{\rm (ix)}]  Formulas for calculating the extremal rank and inertia of
the $A + BC^-B^{*}$, where $C^-$ is a Hermitian generalized inverse of a Hermitian
matrix $C$, \cite{LT-jota,T-ela12}.

\end{enumerate}

Mappings between matrix spaces with symmetric patterns can be constructed arbitrarily,
but the linear function in (\ref{ww18}) is the simplest cases among all matrix maps
with symmetric patterns. The linear matrix inequality in (\ref{ww18}) and its variations are usually
taken as global convex constraints to unknown matrices and vectors
in mathematical programming and optimization theory. Note that the
commonly used definiteness matrices $X \succcurlyeq  0$
 $(X  \preccurlyeq  0)$ is a special case of the inequality in (\ref{ww18}).  Thus, the
 inequality-constraints in (\ref{ww18}) could be regarded as two extensions of definite
 matrix constraints arising in a number of optimization problems
(see, e.g., \cite{Fl,HW,JBAS,TT}). In fact, Problem 1.1 was proposed in the author's
recent paper  \cite{T-mcm}.

The above three problems are closely linked each other. Once analytical formulas for calculating the global maximal
and minimal ranks and inertias in Problem 1.1 are obtained, we can easily use them to solve
Problems 1.2 and 1.3.

The results in the following two lemmas are obvious or well-known (see also \cite{T-laa10,T-laa11} for their
references), which we shall use in the latter part of this paper for solving the previous problems.

\begin{lemma} \label{T11}
Let ${\cal S}$ be a set consisting of  matrices over ${\mathbb
C}^{m\times n},$ and let ${\cal H}$ be a set consisting of Hermitian
matrices over ${\mathbb C}_{{\rm H}}^{m}.$ Then$,$ the following hold$.$
\begin{enumerate}
\item[{\rm (a)}] Under $m =n,$ ${\cal S}$ has a nonsingular matrix if and only if
$\max_{X\in {\cal S}} r(X) = m.$

\item[{\rm (b)}] Under $m =n,$ all $X\in {\cal S}$ are nonsingular if and only if
 $\min_{X\in {\cal S}} r(X) = m.$

\item[{\rm (c)}] $0\in {\cal S}$ if and only if
$\min_{X\in {\cal S}} r(X) = 0.$

\item[{\rm (d)}] ${\cal S} = \{ 0\}$ if and only if
$\max_{X\in {\cal S}} r(X) = 0.$

\item[{\rm (e)}]  All $X\in {\cal S}$ have the same rank if and only if
 $\max_{X\in {\cal S}} r(X) = \min_{X\in {\cal S}} r(X).$

\item[{\rm (f)}] ${\cal H}$ has a matrix $X \succ  0$  $(X \prec  0)$ if and only if
$\max_{X\in {\cal H}} i_{+}(X) = m  \ \left(\max_{X\in
{\cal H}} i_{-}(X) = m \right)\!.$

\item[{\rm (g)}] All $X\in {\cal H}$ satisfy $X \succ 0$ $(X \prec  0),$
namely$,$  ${\cal H}$ is a subset of the cone of positive definite matrices
(negative definite matrices)$,$ if and only if
$\min_{X\in {\cal H}} i_{+}(X) = m \ \left(\min_{X\in {\cal H}}
i_{-}(X) = m\, \right)\!.$

\item[{\rm (h)}] ${\cal H}$ has a matrix  $X \succcurlyeq  0$ $(X \preccurlyeq  0)$ if and only if
$\min_{X \in {\cal H}} i_{-}(X) = 0 \ \left(\min_{X\in
{\cal H}} i_{+}(X) = 0 \,\right)\!.$

\item[{\rm (i)}] All $X\in {\cal H}$ satisfy $X \succcurlyeq  0$ $(X \preccurlyeq  0),$
namely$,$  ${\cal H}$ is a subset of the cone of positive semi-definite matrices (nonpositive definite
matrices)$,$  if and only if $\max_{X\in {\cal H}} i_{-}(X) = 0 \ \left(\max_{X\in {\cal H}}
i_{+}(\,X) = 0\, \right)\!.$

\item[{\rm (j)}]  All $X\in {\cal H}$ have the same positive index of inertia if and only if
 $\max_{X\in {\cal H}} i_{+}(X) = \min_{X\in {\cal H}} i_{+}(X).$

\item[{\rm (k)}]  All $X\in {\cal H}$ have the same negative index of inertia if and only if
 $\max_{X\in {\cal H}} i_{-}(X) = \min_{X\in {\cal H}} i_{-}(X).$
\end{enumerate}
\end{lemma}

The question of whether a given function (matrix map), is
positive or nonnegative (definite or semi-definite) everywhere is
ubiquitous in mathematics and applications. Lemma \ref{T11}(f)--(i)
show that if some explicit formulas for calculating the global
maximal and minimal inertias of a given Hermitian matrix map
are established, we can use them, as demonstrated in Sections 2--5
below, to derive necessary and sufficient conditions for the
Hermitian matrix map to be definite or semi-definite.

\begin{lemma} \label{T12}
Let $A\in {\mathbb C}^m_{{\rm H}},$  $B \in {\mathbb C}^n_{{\rm
H}},$ $Q \in \mathbb C^{m \times n},$ and assume that $P\in \mathbb
C^{m\times m}$ is nonsingular$.$ Then$,$
\begin{eqnarray}
i_{\pm}(PAP^{*}) &\!\! = \!\!&  i_{\pm}(A),
\label{qq17}
\\
 i_{\pm}(\lambda A) &\!\! = \!\!&  \left\{ \begin{array}{ll}  i_{\pm}(A)  &  if  \ \lambda > 0
 \\ i_{\mp}(A)  & if  \ \lambda < 0
 \end{array},
\right.
\label{qq18}
\\
i_{\pm}\!\left[ \begin{array}{cc}  A  & 0 \\  0  & B \end{array}
\right] &\!\! = \!\!&  i_{\pm}(A) + i_{\pm}(B),
\label{qq19}
\\
i_{+}\!\left[ \begin{array}{cc}  0  & Q \\  Q^{*} & 0 \end{array}
\right] &\!\! = \!\!&  i_{-}\!\left[ \begin{array}{cc}  0  & Q \\  Q^{*} & 0 \end{array}
\right] = r(Q).
\label{qq110}
\end{eqnarray}
\end{lemma}

\begin{lemma} [\cite{MS}]\label{T13}
Let $A \in \mathbb C^{m \times n}, \ B \in \mathbb C^{m \times k}$ and
$C \in \mathbb C^{l \times n}.$ Then$,$
the following rank expansion formulas hold
\begin{eqnarray}
r[\, A, \, B \,] &\!\! = \!\!&  r(A) + r(E_AB) = r(B) + r(E_BA),
\label{qq111}
\\
r \!\left[\!\! \begin{array}{c}  A  \\ C  \end{array}
\!\!\right] &\!\! = \!\!&  r(A) + r(CF_A) = r(C) + r(AF_C).
\label{qq112}
\end{eqnarray}
\end{lemma}

\begin{lemma} [\cite{T-laa10}] \label{T14}
Let $A \in {\mathbb C}_{{\rm H}}^{m},$ $B \in \mathbb C^{m\times n},$
 $D \in {\mathbb C}_{{\rm H}}^{n},$   and let
$$
M_1 = \left[\!\!\begin{array}{cc}  A  & B  \\ B^{*}  & 0 \end{array}
\!\!\right]\!, \ \ M_2 = \left[\!\!\begin{array}{cc}  A  & B  \\ B^{*}  & D \end{array}
\!\!\right]\!.
$$
Then$,$ the following expansion formulas hold
\begin{eqnarray}
& i_{\pm}(M_1) = r(B) + i_{\pm}(E_BAE_B), \ \ \ \ \ \ \ \ \ \ \ \ \ \ \ \ \ \ \ \ r(M_1) = 2r(B) + r(E_BAE_B),
\label{qq113}
\\
& i_{\pm}(M_2)  =i_{\pm}(A) + i_{\pm}\!\left[\!\!\begin{array}{cc} 0 &  E_AB
 \\
 B^{*}E_A & D - B^{*}A^{\dag}B \end{array}\!\!\right]\!, \ \
 r(M_2)  = r(A) + r\!\left[\!\!\begin{array}{cc} 0 &  E_AB
 \\
 B^{*}E_A &  D - B^{*}A^{\dag}B  \end{array}\!\!\right]\!.
\label{110}
\end{eqnarray}
In particular$,$ the following hold$.$
\begin{enumerate}
\item[{\rm(a)}] If $A \succcurlyeq  0,$ then
\begin{eqnarray}
i_{+}(M_1)= r[\, A,  \, B \,], \ \ i_{-}(M_1) = r(B), \ \
 r(M_1) = r[\, A,  \, B \,] + r(B).
\label{111j}
\end{eqnarray}

\item[{\rm(b)}] If $A \preccurlyeq  0,$ then
\begin{eqnarray}
i_{+}(M_1) = r(B), \ \ i_{-}(M_1) = r[\, A,  \, B \,], \ \ r(M_1) = r[\,
A, \, B \,] + r(B).
\label{112j}
\end{eqnarray}

\item[{\rm(c)}] If ${\mathscr R}(B) \subseteq {\mathscr R}(A),$ then
\begin{eqnarray}
i_{\pm}(M_2) = i_{\pm}(A) + i_{\pm}(\,  D - B^{*}A^{\dag}B \,), \ \ r(M_2)
= r(A) + r(\,  D - B^{*}A^{\dag}B \,).
\label{113j}
\end{eqnarray}


\item[{\rm (d)}] $i_{\pm}(M_2)  \geqslant i_{\pm}(A) +  i_{\pm}(\, D - B^{*}A^{\dag}B \,)  \geqslant i_{\pm}(A).$

\item[{\rm (e)}] $i_{\pm}(M_1) = m \Leftrightarrow  i_{\mp}(E_BAE_B) =0 \  and  \ r(E_BAE_B) =r(E_B).$

\item[{\rm (f)}]  $i_{+}(M_2) =i_{+}(A) \Leftrightarrow {\mathscr R}(B) \subseteq {\mathscr R}(A) \ and \  D - B^{*}A^{\dag}B \preccurlyeq  0.$

\item[{\rm (g)}]  $i_{-}(M_2) =i_{-}(A) \Leftrightarrow {\mathscr R}(B) \subseteq {\mathscr R}(A) \ and \  D - B^{*}A^{\dag}B \succcurlyeq  0.$

\item[{\rm(h)}] $M_2 \succcurlyeq  0$ $\Leftrightarrow$ $A \succcurlyeq  0,$ ${\mathscr R}(B) \subseteq
{\mathscr R}(A)$ and $D - B^{*}A^{\dag}B \succcurlyeq  0$ $\Leftrightarrow$
$D \succcurlyeq  0,$ ${\mathscr R}(B^{*}) \subseteq {\mathscr R}(D)$ and $ A - BD^{\dag}B^{*} \succcurlyeq  0.$

\item[{\rm(i)}] $M_2 \succ  0$ $\Leftrightarrow A \succ 0$ and $D - B^{*}A^{-1}B \succ 0$
$\Leftrightarrow$ $D \succ 0$  and $A - BD^{-1}B^{*} \succ 0.$
\end{enumerate}
\end{lemma}

Some useful expansion formulas derived from (\ref{qq113}) and (\ref{110})  are
\begin{eqnarray}
 i_{\pm}(\, D - B^{*}A^{\dag}B\,) &\!\! = \!\!& i_{\pm}\!\left[\!\! \begin{array}{cc}  A^{*}AA^{*}  & A^{*}B   \\ B^{*}A^{*} &
 D \end{array} \!\!\right] - i_{\pm}(A),
\label{122}
\\
 i_{\pm}(\, D - B^{*}A^{\dag}B\,) &\!\! = \!\!&  i_{\pm}\!\!\left[\!\! \begin{array}{cc}  A  & B   \\ B^{*} &
 D \end{array} \!\!\right] - i_{\pm}(A) \ \ {\rm if} \ \ {\mathscr R}(B) \subseteq {\mathscr R}(A),
\label{123}
\\
 i_{\pm}\!\left[\!\! \begin{array}{cc} A
 & BF_P \\  F_P B^{*}  & 0 \end{array} \!\!\right] &\! = \!& i_{\pm}\!\left[\!\! \begin{array}{ccc} A
 & B & 0 \\  B^{*} & 0 & P^{*} \\  0 & P & 0 \end{array} \!\!\right] - r(P).
\label{116}
\end{eqnarray}
We shall use them to simplify the inertias of block Hermitian
matrices involving Moore--Penrose inverses of matrices.

\begin{lemma}[\cite{T-lama10}]\label{T15}
Let $A \in {\mathbb C}_{{\rm H}}^{m}$ and $B \in {\mathbb C}^{m\times n}.$ Then$,$
the following expansion formula
$$
\left[\!\!\begin{array}{cc}  A  & B  \\ B^{*}  & 0 \end{array}
\!\!\right]^{\dag} = \left[\!\!\begin{array}{cc}  (E_BAE_B)^{\dag}   & (B^{\dag})^{*} -
(E_BAE_B)^{\dag}E_BA(B^{\dag})^{*}
  \\ B^{\dag} - B^{\dag}AE_B(E_BAE_B)^{\dag}  & - B^{\dag}A(B^{\dag})^{*} +
  B^{\dag}AE_B(E_BAE_B)^{\dag}E_BA(B^{\dag})^{*} \end{array}
\!\!\right]
$$
holds if and only if
$$
r\!\left[\!\!\begin{array}{cc}  A  & B  \\ B^{*}  & 0 \end{array} \!\!\right] =  r[\,A, \, B \,] + r(B),
$$
or
equivalently$,$ $r(E_BAE_B) =r(E_BA).$
\end{lemma}

Solving matrix equations is one of the key problems of matrix computation.
Many techniques were proposed and developed in studying consistency and solutions of
various matrix equations. In this paper, we need the following results on
solvability conditions and general solutions of two simple linear matrix equations.

\begin{lemma} [\cite{KM}] \label{T16}
Let $A, \, B\in \mathbb C^{m\times n}$ be given$.$ Then$,$ the following hold$.$
\begin{enumerate}
\item[{\rm(a)}]  The matrix equation $AX = B$
has a solution $X \in \mathbb C_{{\rm H}}^{n}$ if and only if ${\mathscr R}(B) \subseteq
{\mathscr R}(A)$ and $AB^{*} =BA^{*}.$ In this case$,$ the general
Hermitian solution can be written in the
following parametric form
\begin{equation}
X = A^{\dag}B + (A^{\dag}B)^{*} - A^{\dag}BA^{\dag}A + F_AUF_A,
\label{qq125}
\end{equation}
where $U\in \mathbb C_{{\rm H}}^{n}$ is arbitrary$.$

\item[{\rm(b)}]$AX = B$
has a solution  $0\preccurlyeq  X \in {\mathbb C}_{{\rm H}}^{n}$ if and only if ${\mathscr R}(B) \subseteq {\mathscr
R}(A),$ $AB^{*} \succcurlyeq  0$ and $r(AB^{*}) =r(B).$ In this case$,$ the
general solution $0\preccurlyeq  X \in {\mathbb C}_{{\rm H}}^{n}$  can be written as
\begin{equation}
X = B^{*}(AB^{*})^{\dag}B   + F_AUU^{*}F_A,
 \label{28}
\end{equation}
where $U\in \mathbb C^{n\times n}$ is arbitrary$.$
\end{enumerate}
\end{lemma}

\begin{lemma} \label{T17}
Let $A \in \mathbb C^{m \times n}$ and $B \in \mathbb C_{{\rm
H}}^{m}$ be given$.$ Then$,$ the following hold$.$
\begin{enumerate}
\item[{\rm(a)}] {\rm \cite{Gr}} The matrix equation
\begin{equation}
AXA^{*} = B
\label{29}
\end{equation}
has a solution $X \in \mathbb C_{{\rm H}}^{n}$ if and only if
${\mathscr R}(B) \subseteq {\mathscr R}(A),$ or equivalently$,$
$AA^{\dag}B = B.$

\item[{\rm(b)}] {\rm \cite{T-laa10}} Under ${\mathscr R}(B) \subseteq {\mathscr R}(A),$
the general Hermitian solution of {\rm (\ref{29})} can be written in the following two forms
\begin{eqnarray}
X  = A^{\dag}B(A^{\dag})^{*} +  F_AV + V^{*}F_A,
\label{z43}
\end{eqnarray}
respectively$,$ where $V \in {\mathbb C}^{n \times n}$ is arbitrary$.$

\item[{\rm(c)}] {\rm \cite{Gr,KM}} The matrix equation
\begin{equation}
AXA^{*} = B
\label{211}
\end{equation}
has a solution $0\preccurlyeq  X \in {\mathbb C}_{{\rm H}}^{n}$  if and only if
$B \succcurlyeq  0$  and ${\mathscr R}(B) \subseteq {\mathscr R}(A).$ In this case$,$
the general positive semi-definite solution of {\rm (\ref{211})} can be written in
the following parametric form
\begin{equation}
X = A^{-}B(A^{-})^{*} +  F_AUU^{*}F_A = (\, A^{\dag} + F_AV\,)B (\, A^{\dag} + F_AV\,)^{*} +  F_AUU^{*}F_A,
\label{212}
\end{equation}
where $A^-$ is an arbitrary g-inverse of $A,$ and $V\in {\mathbb C}^{n \times m}$ and
$U\in \mathbb C^{n \times n}$ are arbitrary$.$
\end{enumerate}
\end{lemma}

In order to simplify various matrix-valued function involving generalized inverse of matrices and
arbitrary matrices, we need the following results on ranks of matrices.

\begin{lemma}[\cite{LT-jota,T-laa11}]\label{T18}
Let $A \in {\mathbb C}_{{\rm H}}^{m}$, $B \in {\mathbb C}^{m\times n}$ and $C \in {\mathbb C}^{p \times m}$ be
 given$.$ Then$,$   the global maximal and minimal ranks and inertias of $A-BXC-(BXC)^{*}$ are given by
\begin{align}
& \max_{X\in {\mathbb C}^{n \times p}}\!\!\!r[\, A-BXC-(BXC)^{*}\,] =  \min
\left\{ r[\,A,\, B,\, C^{*}\,],
 \ \ r\!\left[\begin{array}{cc} A & B
\\B^{*} & 0
\end{array}\right]\!, \ \ r\!\left[\!\begin{array}{cc} A & C^{*}
\\ C & 0
\end{array}\!\right] \right\}\!,
\label{129}
\\
& \min_{X\in {\mathbb C}^{n \times p}}\!\!\!r[\, A-BXC-(BXC)^{*}\,] =
2r[\,A,\, B,\, C^{*}\,] + \max\{\, s_{+} + s_{-},  \ t_{+} + t_{-},
\ s_{+} + t_{-}, \ s_{-} + t_{+} \, \}, \label{130}
\\
& \max_{X\in {\mathbb C}^{n \times p}}\!\!\!i_{\pm}[\, A-BXC-(BXC)^{*}\,] = \min\!\left\{i_{\pm}\!\left[\!\begin{array}{ccc} A
 & B  \\  B^{*}  & 0
   \end{array}\!\right], \ \ i_{\pm}\!\left[\!\begin{array}{ccc} A
 & C^{*}  \\  C  & 0
   \end{array}\!\right] \right\}\!,
\label{131}
\\
& \min_{X\in {\mathbb C}^{n \times p}}\!\!\!i_{\pm}[\, A-BXC-(BXC)^{*}\,] =
r[\,A,\, B,\, C^{*}\,] + \max\{\, s_{\pm}, \ \ t_{\pm} \, \},
\label{132}
\end{align}
where
$$
s_{\pm}  = i_{\pm}\!\left[\!\!\begin{array}{cc} A & B \\ B^{*} & 0\end{array}\!\!\right]
 - r\!\left[\begin{array}{ccc} A & B  & C^{*} \\ B^{*} & 0 & 0
\end{array}\!\!\right]\!, \ \
 t_{\pm}  =i_{\pm}\!\left[\!\!\begin{array}{cc} A & C^{*} \\ C & 0\end{array}\!\!\right]
 - r\!\left[\begin{array}{ccc} A & B  & C^{*} \\ C & 0 & 0
\end{array}\!\!\right]\!.
$$
In particular, if ${\mathscr R}(B)\subseteq {\mathscr R}(C^{*}),$ then
\begin{align}
& \max_{X\in {\mathbb C}^{n \times p}}\!\!\!r[\, A-BXC-(BXC)^{*}\,] =  \min
\left\{ r[\,A,\, C^{*}\,],
 \ \ r\!\left[\begin{array}{cc} A & B
\\B^{*} & 0
\end{array}\right] \right\}\!,
\label{133}
\\
& \min_{X\in {\mathbb C}^{n \times p}}\!\!\!r[\, A-BXC-(BXC)^{*}\,] =
2r[\,A,\, C^{*}\,] + r\!\left[\begin{array}{cc} A & B
\\B^{*} & 0
\end{array}\right] - 2r\!\left[\begin{array}{cc} A & B
\\C & 0
\end{array}\right]\!,
\label{134}
\\
& \max_{X\in {\mathbb C}^{n \times p}}\!\!\!i_{\pm}[\, A-BXC-(BXC)^{*}\,] = i_{\pm}\!\left[\!\begin{array}{ccc} A
 & B  \\  B^{*}  & 0   \end{array}\!\right]\!,
\label{135}
\\
& \min_{X\in {\mathbb C}^{n \times p}}\!\!\!i_{\pm}[\, A-BXC-(BXC)^{*}\,] =
r[\,A,\, C^{*}\,] + i_{\pm}\!\left[\!\begin{array}{ccc} A
 & B  \\  B^{*}  & 0
   \end{array}\!\right]- r\!\left[\begin{array}{cc} A & B
\\ C & 0
\end{array}\right]\!,
\label{136}
\end{align}
and
\begin{align}
& \max_{X\in {\mathbb C}^{n \times m}}\!\!\!r[\,A - BX - (BX)^{*}\,] = \min \left\{m,
 \ \ r\!\left[\begin{array}{cc} A & B
\\B^{*} & 0
\end{array}\right] \right\}\!,
\label{137}
\\
& \min_{X\in {\mathbb C}^{n \times m}}\!\!\!r[\,A - BX - (BX)^{*}\,] =  r\!\left[\begin{array}{cc} A & B
\\B^{*} & 0
\end{array}\right]-2r(B),
\label{138}
\\
& \max_{X\in {\mathbb C}^{n \times m}}\!\!\!i_{\pm}[\,A  - BX - (BX)^{*}\,] =
i_{\pm}\!\left[\!\begin{array}{ccc} A & B  \\  B^{*}  & 0
   \end{array}\!\right]\!,
\label{139}
\\
& \min_{X\in {\mathbb C}^{n \times m}}\!\!\!i_{\pm}[\,A -  BX - (BX)^{*}\,] =  i_{\pm}\!\left[\!\begin{array}{ccc} A
 & B  \\  B^{*}  & 0
   \end{array}\!\right]- r(B).
\label{140}
\end{align}
\end{lemma}

This paper is organized as follows. In Section 2, we derive general Hermitian solution of
the LMI  in (\ref{ww18}) by using generalized inverses of matrices, Lemmas 1.9 and
1.10, and present some algebraic properties of the Hermitian solutions. In Sections 3--5,
 we derive explicit solutions to Problems 1.1--1.3, and present various consequences of the rank
 and inertia formulas obtained. In Section 6, we calculate the global maximal and minimal ranks
 and inertias of the Hermitian solution of $BXB^{*} \succcurlyeq  A$, as well as the global maximal and minimal ranks
and inertias of the submatrices in a Hermitian solution of $BXB^{*} \succcurlyeq  A$.

\section{General Hermitian solutions of the LMIs $BXB^{*} \succcurlyeq  (\succ , \,  \preccurlyeq , \,  \prec ) \, A$ and their
properties}
\renewcommand{\theequation}{\thesection.\arabic{equation}}
\setcounter{section}{2}
\setcounter{equation}{0}

Concerning the global maximal and minimal  ranks and inertias of  (\ref{ww18}),
we have the following known result.

\begin{lemma} [\cite{T-laa10,TL}] \label{T21}
Let  $A \in {\mathbb C}_{{\rm H}}^{m}$ and $B \in \mathbb C^{m\times n}$ be given$,$ and  define
$M = \!\left[\!\! \begin{array}{cc} A & B \\ B^{*} & 0 \end{array} \!\!\right].$
Then$,$
 the global maximal and minimal rank and inertias of $A-BXB^{*}$ are given by
\begin{align}
&  \max_{X \in {\mathbb C}_{{\rm H}}^{n}}\!r(\, A - BXB^{*} \,)  = r[\,A, \, B\,],
\\
& \min_{X \in {\mathbb C}_{{\rm H}}^{n}} \!r(\, A - BXB^{*} \,) =
2r[\, A, \, B \,]  - r(M),
\label{qq21}
\\
 & \max_{X\in {\mathbb C}_{{\rm H}}^{n}}\!i_{\pm}(\, A - BXB^{*} \,)
 = i_{\pm}(M),
\\
& \min_{X \in {\mathbb C}_{{\rm H}}^{n}} \! i_{\pm}(\, A - BXB^{*} \,)  =  r[\, A, \, B \,]  -
 i_{\mp}(M).
 \label{qq22}
\end{align}
\end{lemma}


We next solve the two inequalities  in (\ref{ww17}) and give their general Hermitian solutions by
using Lemmas 1.9 and 1.10, some partial conclusions were given in \cite{T-mia}.

\begin{theorem} \label{T22}
Let $A \in {\mathbb C}_{{\rm H}}^{m}$ and $B \in \mathbb C^{m\times n}$ be given$.$  Then$,$ the following hold$.$
\begin{enumerate}
\item[{\rm(a)}] There exists an $X \in {\mathbb C}^{n}_{{\rm H}}$ such that
\begin{equation}
BXB^{*}  \succcurlyeq  A
\label{qq26}
\end{equation}
if and only if
\begin{equation}
E_BAE_B \preccurlyeq  0 \ \ \ and  \ \ \  r(E_BAE_B)= r(E_BA),
\label{qq27}
\end{equation}
or equivalently$,$
\begin{equation}
i_{+}(M) = r(B)  \ \ \ and  \ \ \  i_{-}(M) = r[\, A, \, B\,].
\label{qq27a}
\end{equation}
In this case$,$ the general Hermitian solution of {\rm (\ref{qq26})} can be written as
\begin{align}
& X  = B^{\dag}A(B^{\dag})^{*} - B^{\dag}AE_B(E_BAE_B)^{\dag}E_BA(B^{\dag})^{*} + UU^{*}
 + F_BV + V^{*}F_B,
\label{qq28}
\\
& BXB^{*}  =  A - AE_B(E_BAE_B)^{\dag}E_BA + BUU^{*}B^{*},
\label{qq28a}
\\
& A - BXB^{*} = AE_B(E_BAE_B)^{\dag}E_BA - BUU^{*}B^{*},
\label{qq28ab}
\end{align}
where $U, \, V\in {\mathbb C}^{n \times n}$ are arbitrary$.$

\item[{\rm(b)}] There exists an $X \in {\mathbb C}^{n}_{{\rm H}}$ such that
\begin{equation}
BXB^{*} \succ  A
\label{qq29}
\end{equation}
if and only if
\begin{equation}
E_BAE_B \preccurlyeq  0 \ \ \ and  \ \  r(E_BAE_B)= r(E_B)
\label{qq210}
\end{equation}
hold$.$ In this case$,$ the general Hermitian solution of {\rm (\ref{qq29})} can be
written as {\rm (\ref{qq28})}$,$ in which  $U$ is a matrix such that
$r[\, AE_B(E_BAE_B)^{\dag}E_BA -  BUU^{*}B^{*}\,] =m,$ and $V\in {\mathbb C}^{n \times n}$
is arbitrary$.$
\end{enumerate}
In particular$,$ the following hold$.$
\begin{enumerate}
\item[{\rm(c)}] If $BXB^{*} = A$ is consistent$,$ then the general
Hermitian solution of $BXB^{*}  \succcurlyeq  A$ can be written as
\begin{equation}
X = B^{\dag}A(B^{\dag})^{*} + UU^{*} + F_BV + V^{*}F_B,
\label{qq216}
\end{equation}
where $U, \, V \in \mathbb C^{n \times n}$ are arbitrary$.$

\item[{\rm(d)}] If $BXB^{*} = A$ is consistent$,$ then  $BXB^{*} \succ  A$ has a Hermitian solution if and only if
$r(B) = m,$ in which case$,$ the general Hermitian solution of  the LMI
 can be written as
\begin{equation}
X = B^{\dag}A(B^{\dag})^{*}  + U + F_BV + V^{*}F_B,
\label{qq217}
\end{equation}
where  $U \in {\mathbb C}_{{\rm H}}^{n}$ is
arbitrary matrix such that $BUB^{*} \succ 0,$  and $V \in \mathbb C^{n \times n}$  is arbitrary$.$
\end{enumerate}
\end{theorem}

\noindent {\bf Proof.} \ It is obvious that (\ref{qq26}) is equivalent to
\begin{equation}
BXB^{*} = A + YY^{*}
\label{qq220}
\end{equation}
for some matrix $Y$.  In other words, (\ref{qq26}) can be relaxed to a matrix equation with two
unknown matrices. We obtain from Lemma 1.10(a) that (\ref{qq220}) is solvable for
$X \in {\mathbb C}^{n}_{{\rm H}}$ if and only if $E_B(A + YY^{*}) = 0,$ that is,
\begin{equation}
E_BYY^{*} = -E_BA.
\label{qq221}
\end{equation}
From Lemma 1.9(b), (\ref{qq221}) is solvable for $YY^{*}$ if and only if $E_BAE_B \preccurlyeq  0$
and $r(E_BAE_B)= r(E_BA)$, establishing (\ref{qq27}), which is further equivalent to (\ref{qq27a})
by (\ref{qq111}) and (\ref{qq113}). In this case, the general positive semi-definite solution of
(\ref{qq221}) can be written as
\begin{equation}
YY^{*} =  -AE_B(E_BAE_B)^{\dag}E_BA +  BB^{\dag}UU^{*}BB^{\dag},
\label{qq222}
\end{equation}
where $U \in {\mathbb C}^{m \times m}$ is arbitrary. Substituting
the $YY^{*}$ into (\ref{qq220}) gives
\begin{equation}
BXB^{*} = A -AE_B(E_BAE_B)^{\dag}E_BA +  BB^{\dag}UU^{*}BB^{\dag}.
\label{qq223}
\end{equation}
By Lemma 1.10(b), the general Hermitian solution of (\ref{qq223}) can be written as
\begin{equation}
X = B^{\dag}A(B^{\dag})^{*} - B^{\dag}AE_B(E_BAE_B)^{\dag}E_BA(B^{\dag})^{*} + B^{\dag}UU^{*}(B^{\dag})^{*}
 + F_BV + V^{*}F_B,
\label{qq224}
\end{equation}
where $V \in \mathbb C^{n \times n}$ is arbitrary$.$ Replacing the matrix
$B^{\dag}UU^{*}(B^{\dag})^{*}$ in (\ref{qq224}) with $UU^{*}$ yields (\ref{qq28}), which
is also the general Hermitian solution of (\ref{qq26}).

Substituting (\ref{qq28}) into $A - BXB^{*}$ gives
\begin{eqnarray}
A -BXB^{*} &\!\! = \!\!&  A -  BB^{\dag}ABB^{\dag} + BB^{\dag}AE_B(E_BAE_B)^{\dag}E_BABB^{\dag} - BUU^{*}B^{*}  \nb
\\
&\!\! = \!\!&  AE_B(E_BAE_B)^{\dag}E_BA  - BUU^{*}B^{*},
\label{qq224a}
\end{eqnarray}
Note that $AE_B(E_BAE_B)^{\dag}E_BA \preccurlyeq  0$. Then, we have
\begin{eqnarray*}
& i_{-}[\, AE_B(E_BAE_B)^{\dag}E_BA -  BUUB^{*} \,] =
r[\, AE_B(E_BAE_B)^{\dag}E_BA -  BUUB^{*} \,]
\\
& =  r[\, AE_B(E_BAE_B)^{\dag}E_BA, \,  BUUB^{*}\,] = r[\, AE_B, \, BU\,].
\end{eqnarray*}
In consequence,
$$
\max_{U\succcurlyeq  0} i_{-}[\, AE_B(E_BAE_B)^{\dag}E_BA -  BUU^{*}B^{*} \,] =  \max_{U}r[\, AE_B, \,   BU\,]
 = r[\, AE_B, \,  B\,] = r(E_BAE_B) + r(B),
$$
so that (\ref{qq29}) holds if and only if $r(E_BAE_B) + r(B) = m,$
establishing (b). If $BXB^{*} = A$ is consistent$,$ then $E_B A = 0$.
In this case, (a) and (b) reduce to (c) and (d).  \qquad $\Box$

\medskip

Replacing $A$ with $-A$, and $X$ with $-X$ in Theorem \ref{T22} leads to
 the following consequence.

\begin{corollary} \label{T22a}
Let $A \in {\mathbb C}_{{\rm H}}^{m}$ and $B \in \mathbb C^{m\times n}$ be given$.$
Then$,$ the following hold$.$
\begin{enumerate}
\item[{\rm(a)}] There exists an $X \in {\mathbb C}^{n}_{{\rm H}}$ such that
\begin{equation}
BXB^{*}  \preccurlyeq  A
\label{qq211}
\end{equation}
if and only if
\begin{equation}
E_BAE_B \succcurlyeq  0 \ \ \ and  \ \ \  r(E_BAE_B)= r(E_BA),
\label{qq212}
\end{equation}
or equivalently$,$
\begin{equation}
i_{+}(M) = r[\, A, \, B\,] \ \ \ and  \ \ \ i_{-}(M) = r(B).
\label{qq212a}
\end{equation}
 In this case$,$ the general Hermitian solution of {\rm (\ref{qq211})}
can be written in the following parametric form
\begin{align}
& X =  B^{\dag}A(B^{\dag})^{*} - B^{\dag}AE_B(E_BAE_B)^{\dag}E_BA(B^{\dag})^{*} - UU^{*}
 + F_BV + V^{*}F_B,
\label{qq212b}
 \\
&  BXB^{*}  = A - AE_B(E_BAE_B)^{\dag}E_BA - BUU^{*}B^{*},
\label{qq212c}
\\
&  A - BXB^{*}  = AE_B(E_BAE_B)^{\dag}E_BA  + BUU^{*}B^{*},
\label{qq212d}
\end{align}
where  $U, \, V \in \mathbb C^{n \times n}$ are arbitrary$.$

\item[{\rm(b)}] There exists an $X \in {\mathbb C}_{{\rm H}}^{n}$ such that
\begin{equation}
BXB^{*} \prec  A
 \label{qq214}
\end{equation}
if and only if
\begin{equation}
E_BAE_B \succcurlyeq  0 \ \ \ and  \ \  r(E_BAE_B)= r(E_B)
\label{qq215}
\end{equation}
hold$.$  In this case$,$ the general Hermitian solution of
{\rm (\ref{qq214})} can be written as {\rm (\ref{qq212b})}$,$ in which
$U$ is a matrix such that $r[\, AE_B(E_BAE_B)^{\dag}E_BA -  BUU^{*}B^{*}\,] =m,$
and $V \in \mathbb C^{n \times n}$ is arbitrary$.$
\end{enumerate}
In particular$,$
\begin{enumerate}
\item[{\rm(c)}] If $BXB^{*} = A$ is consistent$,$ then the general Hermitian solution of
$BXB^{*}  \preccurlyeq  A$ can be written as
\begin{equation}
X = B^{\dag}A(B^{\dag})^{*} - UU^{*} + F_BV + V^{*}F_B,
\label{qq218}
\end{equation}
where  $U, \, V \in \mathbb C^{n \times n}$ are arbitrary$.$

\item[{\rm(d)}] If $BXB^{*} = A$ is consistent$,$ then of $BXB^{*} \prec  A$  has a Hermitian solution if
and only if $r(B) = m,$ in which case$,$ the general Hermitian solution of the LMI
 can be written as
\begin{equation}
X = B^{\dag}A(B^{\dag})^{*} -  U + F_BV + V^{*}F_B,
\label{qq219}
\end{equation}
where  $U \in {\mathbb C}_{{\rm H}}^{n}$ is
arbitrary matrix such that $BUB^{*} \succ 0,$ and $V \in \mathbb C^{n \times n}$  is arbitrary$.$
\end{enumerate}
\end{corollary}

We next establish some algebraic properties of the fixed part in (\ref{qq28}) and (\ref{qq212b}).

\begin{corollary} \label{T23}
Let $A \in {\mathbb C}_{{\rm H}}^{m}$ and  $B \in \mathbb C^{m\times n}$ be given$,$ and let
\begin{equation}
\widehat{X} = B^{\dag}A(B^{\dag})^{*} - B^{\dag}AE_B(E_BAE_B)^{\dag}E_BA(B^{\dag})^{*}.
\label{qq226}
\end{equation}
Then$,$ the following hold$.$
\begin{enumerate}
\item[{\rm(a)}] Under the condition in {\rm (\ref{qq27})}$,$
\begin{enumerate}
\item[{\rm(i)}] $\widehat{X}$ is a Hermitian solution of $BXB^{*} \succcurlyeq  A.$

\item[{\rm(ii)}] $\widehat{X}$ can be written as $\widehat{X} = [\, 0, \, I_n \,]\left[\begin{array}{cc} - A & B
\\B^{*} & 0\end{array}\right]^{\dag}\left[\begin{array}{cc} 0 \\ I_n
\end{array}\right]\!.$

\item[{\rm(iii)}] $\widehat{X}$ satisfies the following equalities
\begin{eqnarray}
& &  i_{+}(\widehat{X}) = i_{+}(B\widehat{X}B^{*}) = i_{+}(A),
\label{230jj}
\\
& & i_{-}(\widehat{X}) =  i_{-}(B\widehat{X}B^{*}) = i_{-}(A) + r(B) - r[\, A, \, B \,],
\label{231jj}
\\
& & r(\widehat{X}) =  r(B\widehat{X}B^{*}) = r(A) + r(B) - r[\, A, \, B \,],
\label{232jj}
\\
& & i_{-}(\,A - B\widehat{X}B^{*}\,) = r(\,A - B\widehat{X}B^{*}\,) =
 r(A) - r(\,B\widehat{X}B^{*}\,) = r[\, A, \, B\,] - r(B),
\label{233jj}
\\
& & \max_{\succcurlyeq }\{ A - BXB^{*} \ | \  BXB^{*}  \succcurlyeq  A \  and \
 X \in {\mathbb C}_{{\rm H}}^{n} \}  = A - B\widehat{X}B^{*}.
\label{234jj}
\end{eqnarray}
\end{enumerate}

\item[{\rm(b)}] Under the condition in {\rm (\ref{qq212})}$,$
\begin{enumerate}
\item[{\rm(i)}] $\widehat{X}$ is a Hermitian solution of $BXB^{*} \preccurlyeq  A.$

\item[{\rm(ii)}] $\widehat{X}$ can be written as $\widehat{X} = [\, 0, \, I_n \,]\left[\begin{array}{cc} - A & B
\\B^{*} & 0\end{array}\right]^{\dag}\left[\begin{array}{cc} 0 \\ I_n
\end{array}\right]\!.$

\item[{\rm(iii)}] $\widehat{X}$ satisfies the following equalities
\begin{align}
&  i_{+}(\widehat{X}) =  i_{+}(B\widehat{X}B^{*}) =
 i_{+}(A) + r(B) - r[\, A, \, B \,],
\label{235jj}
\\
&  i_{-}(\widehat{X}) =  i_{-}(B\widehat{X}B^{*}) = i_{-}(A),
\label{236jj}
\\
&  r(\widehat{X}) = r(B\widehat{X}B^{*}) = r(A) + r(B) - r[\, A, \, B \,],
\label{237jj}
\\
&  i_{+}(\,A - B\widehat{X}B^{*}\,) = r(\,A - B\widehat{X}B^{*}\,)
= r(A) - r(\,B\widehat{X}B^{*}\,) = r[\, A, \, B\,] - r(B),
\label{238jj}
\\
&  \min_{\succcurlyeq }\{ A - BXB^{*} \ | \  BXB^{*}  \preccurlyeq  A \  and \
 X \in {\mathbb C}_{{\rm H}}^{n} \}  = A - B\widehat{X}B^{*}.
 \label{239jj}
\end{align}
\end{enumerate}
\end{enumerate}
\end{corollary}

\noindent {\bf Proof.} \ Under the condition in (\ref{qq26}), comparing (\ref{qq226}) with Lemma 1.5
leads to (ii) of (a).

Applying (\ref{123}) to (\ref{qq226}) and simplifying  by congruence matrix operations,
 we obtain
\begin{align}
 i_{\pm}(\widehat{X}) &=  i_{\pm}[\,  B^{\dag}A(B^{\dag})^{*} - B^{\dag}AE_B(E_BAE_B)^{\dag}E_BA(B^{\dag})^{*} \,] \nb
\\
&  =  i_{\pm}\!\left[\!\begin{array}{cc} E_BAE_B &  E_BA(B^{\dag})^{*} \\
B^{\dag}AE_B & B^{\dag}A(B^{\dag})^{*} \end{array}\!\right] - i_{\pm}(E_BAE_B) \nb
\\
&  =  i_{\pm}\!\left[\!\begin{array}{cc} A &  A(B^{\dag})^{*} \\
B^{\dag}A & B^{\dag}A(B^{\dag})^{*} \end{array}\!\right] - i_{\pm}(E_BAE_B)  = i_{\pm}(A) - i_{\pm}(E_BAE_B).
\label{qq227}
\end{align}
In consequence,
$$
i_{+}(\widehat{X}) = i_{+}(A), \ \ i_{-}(\widehat{X}) = i_{-}(A) -
i_{-}(E_BAE_B) = i_{-}(A) -  r(E_BA) = i_{-}(A) + r(B) - r[\, A, \, B \,],
$$
establishing (iii) of (a).

Under the condition in (\ref{qq27}), applying (\ref{123}) and simplifying  by congruence matrix operations,
we obtain
\begin{align}
i_{\pm}(\,A - B\widehat{X}B^{*}\,) & = i_{\pm}[\,  A - BB^{\dag}ABB^{\dag} +
BB^{\dag}AE_B(E_BAE_B)^{\dag}E_BABB^{\dag} \,] \nb
\\
&  =  i_{\pm}\!\left[\!\begin{array}{cc} - E_BAE_B &  E_BABB^{\dag} \\
BB^{\dag}AE_B &  A - BB^{\dag}ABB^{\dag} \end{array}\!\right] - i_{\mp}(E_BAE_B) \nb
\\
& =  i_{\pm}\!\left[\!\begin{array}{cc} - E_BAE_B &  E_BA \\
AE_B & 0 \end{array}\!\right] - i_{\mp}(E_BAE_B) \nb
\\
& =  r(E_BA) - i_{\mp}(E_BAE_B).
\label{qq229}
\end{align}
In consequence,
\begin{align*}
& i_{+}(\,A - B\widehat{X}B^{*}\,) = r(E_BA) - i_{-}(E_BAE_B) =r(E_BA) - r(E_BAE_B) = 0,
\\
& i_{-}(\,A - B\widehat{X}B^{*}\,) = r(E_BA) - i_{+}(E_BAE_B) =r[\, A, \, B\,] - r(B),
\end{align*}
establishing (iv) of (a).

Substituting (\ref{qq28}) into $A -BXB^{*}$ gives
\begin{equation}
A - BXB^{*} = A - B\widehat{X}B^{*} - BUU^{*}B^{*} \preccurlyeq   A -
B\widehat{X}B^{*} \label{qq228}
\end{equation}
for any $U\in \mathbb C^{n\times n}$, which implies (v) of (a).
 Result (b) can be shown similarly. \qquad $\Box$

\begin{corollary} \label{TS25}
Let $A \in {\mathbb C}_{{\rm H}}^{m}$ and $B \in \mathbb C^{m\times n}$ be given$.$
\begin{enumerate}
\item[{\rm(a)}] Assume that {\rm (\ref{qq26})} has a solution$,$ and define
\begin{equation}
{\cal S}_1 =\{\, X\in {\mathbb C}_{{\rm H}}^{n} \ | \  BXB^{*} \succcurlyeq   A \, \}.
\label{427yy}
\end{equation}
Then$,$ the minimal matrices of $BXB^{*}$ and $BXB^{*} - A$ subject to $X \in {\cal S}_1$
in the L\"owner partial ordering are given by
\begin{eqnarray}
& &\min_{\succcurlyeq }\{ BXB^{*} \ | \  X \in {\cal S}_1 \}  =  A - AE_B(E_BAE_B)^{\dag}E_BA,
\label{428yy}
\\
& & \min_{\succcurlyeq }\{ BXB^{*} - A \ | \  X \in {\cal S}_1 \}  =  - AE_B(E_BAE_B)^{\dag}E_BA,
\label{429yy}
\end{eqnarray}
while the extremal ranks and inertias of $BXB^{*}$ and $BXB^{*} -A$ subject
  to $X \in {\cal S}_1$ are given by
\begin{align}
& \max_{X \in {\cal S}_1} \!r(BXB^{*})  = \max_{X \in {\cal S}_1} i_{+}(BXB^{*})  = r(B),
\label{430yy}
\\
& \min_{X \in {\cal S}_1}\! r(BXB^{*}) = \min_{X \in {\cal S}_1} i_{+}(BXB^{*}) = i_{+}(A),
\label{431yy}
\\
& \max_{X \in {\cal S}_1} \!i_{-}(BXB^{*})  = r(B) + i_{-}(A) - r[\, A, \, B \,],
\label{432yy}
\\
& \min_{X \in {\cal S}_1} \!i_{-}(BXB^{*})   =  0,
\label{433yy}
\\
& \max_{X \in {\cal S}_1} \!r(\,BXB^{*} - A\,)   = r[\, A, \, B \,],
\label{434yy}
\\
&  \min_{X \in {\cal S}_1} \!r(\, BXB^{*} - A  \,)  =  r[\, A, \, B \,] - r(B).
\label{435yy}
\end{align}

\item[{\rm(b)}]  Assume that {\rm (\ref{qq211})} has a solution$,$ and define
\begin{equation}
{\cal S}_2 =\{\, X\in {\mathbb C}_{{\rm H}}^{n} \ | \  BXB^{*} \preccurlyeq   A \, \}.
\label{436yy}
\end{equation}
Then$,$ the maximal matrices of $BXB^{*}$ and $BXB^{*} - A$ subject to $X \in {\cal S}_2$
in the L\"owner partial ordering are given by
\begin{align}
& \max_{\succcurlyeq }\{ BXB^{*} \ | \  X \in {\cal S}_2 \}  =  A - AE_B(E_BAE_B)^{\dag}E_BA,
\label{437yy}
\\
& \max_{\succcurlyeq }\{ BXB^{*} - A \ | \  X \in {\cal S}_2 \}  =  - AE_B(E_BAE_B)^{\dag}E_BA,
\label{438yy}
\end{align}
while
the extremal ranks and inertias of $BXB^{*}$ and  $BXB^{*} -A$ subject to $X \in {\cal S}_2$ are given by
\begin{align}
&\max_{X \in {\cal S}_2} \!r(BXB^{*}) = \max_{X \in {\cal S}_2} i_{-}(BXB^{*})  = r(B),
\label{439yy}
\\
& \min_{X \in {\cal S}_2}\! r(BXB^{*}) = \min_{X \in {\cal S}_2} i_{-}(BXB^{*}) = i_{-}(A),
\label{440yy}
\\
& \max_{X \in {\cal S}_2} \!i_{+}(BXB^{*}) =  r(B) + i_{+}(A) - r[\, A, \, B \,],
\label{441yy}
\\
& \min_{X \in {\cal S}_2} \!i_{+}(BXB^{*})  =  0,
\label{442yy}
\\
& \max_{X \in {\cal S}_2} \!r(\, A - BXB^{*} \,) = r[\, A, \, B \,],
\label{443yy}
\\
& \min_{X \in {\cal S}_2} \!r(\, A - BXB^{*} \,)  =  r[\, A, \, B \,] - r(B).
\label{444yy}
\end{align}
\end{enumerate}
\end{corollary}

\noindent {\bf Proof.} \   It can be seen from (\ref{qq28a}) that
\begin{align}
 BXB^{*} \succcurlyeq  A  - AE_B(E_BAE_B)^{\dag}E_BA, \ \ \ BXB^{*} - A \succcurlyeq   - AE_B(E_BAE_B)^{\dag}E_BA
\label{447yy}
\end{align}
hold for any $U\in {\mathbb C}^{n\times n}$, which  implies (\ref{428yy}) and  (\ref{429yy}).
Applying (\ref{256})--(\ref{260}) to (\ref{qq28a}) and
simplifying  by congruence matrix operations, we obtain
 \begin{align*}
 \max_{X \in {\cal S}_1} \!r(BXB^{*}) & =  \max_{U\in {\mathbb C}^{n\times n}} r(\, B\widehat{X}B^{*} +  BUU^{*}B^{*} \,)
 = r[\, B, \, B\widehat{X}B^{*} \,] = r(B),
\\
\min_{X \in {\cal S}_1} \!r(BXB^{*}) & = \min_{U\in {\mathbb C}^{n\times n}} r(\, B\widehat{X}B^{*} +  BUU^{*}B^{*} \,) \nb
\\
&=  i_{+}(B\widehat{X}B^{*}) +  r[\, B, \, B\widehat{X}B^{*} \,] - i_{+}\!\left[\!\! \begin{array}{cc}
B\widehat{X}B^{*} & B \\  B^{*} & 0 \end{array}
\!\!\right]  = i_{+}(B\widehat{X}B^{*}) = i_{+}(A),
\\
 \max_{X \in {\cal S}_1}\! i_{+}(BXB^{*}) & =  \max_{U\in {\mathbb C}^{n\times n}} i_{+}(\, B\widehat{X}B^{*} +  BUU^{*}B^{*} \,)
 =  i_{+}\!\left[\!\! \begin{array}{cc} B\widehat{X}B^{*} & B \\  B^{*} & 0 \end{array}
\!\!\right] = r(B),
\\
 \min_{X \in {\cal S}_1} \!i_{+}(BXB^{*}) & = \min_{U\in {\mathbb C}^{n\times n}} i_{+}(\, B\widehat{X}B^{*} +  BUU^{*}B^{*} \,)
  = i_{+}(B\widehat{X}B^{*}) = i_{+}(A),
\\
 \max_{X \in {\cal S}_1} \!i_{-}(BXB^{*}) &  =
 \max_{U\in {\mathbb C}^{n\times n}} i_{-}(\, B\widehat{X}B^{*} +  BUU^{*}B^{*} \,)  =  i_{-}(B\widehat{X}B^{*})
 = r(B) + i_{-}(A) - r[\, B, \, A \,],
\\
 \min_{X \in {\cal S}_1}\! i_{-}(BXB^{*}) &  =
 \min_{U\in {\mathbb C}^{n\times n}} i_{-}(\, B\widehat{X}B^{*} +  BUU^{*}B^{*} \,)  = r[\, B, \, B\widehat{X}B^{*}\,] -
  i_{+}\!\left[\!\! \begin{array}{cc} B\widehat{X}B^{*} & B \\  B^{*}  & 0 \end{array}
\!\!\right]  = 0,
\end{align*}
establishing (\ref{428yy})--(\ref{431yy}). Applying (\ref{qq111}),
(\ref{qq112}) to (\ref{qq28a}),
we obtain
\begin{align*}
& \max_{X \in {\cal S}_1} \!r(\,  A - BXB^{*} \,) =
\max_{U\in {\mathbb C}^{n\times n}} r[\, AE_B, \, BU \,] = r[\,  AE_B, \, B \,] =
r\!\left[\!\! \begin{array}{cc} A & B \\  B^{*}  & 0 \end{array}
\!\!\right] - r(B) =  r[\, A, \, B \,],
\\
& \min_{X \in {\cal S}_1} \!r(\,A -  BXB^{*} \,) =
\min_{U\in {\mathbb C}^{n\times n}} r[\,  AE_B, \, BU \,] = r(AE_B) =  r[\, A, \, B \,] - r(B),
\end{align*}
establishing (\ref{434yy}) and (\ref{435yy}). Result (b) can be shown similarly. \qquad $\Box$

\medskip

Some direct consequences of Theorem \ref{T22} are given below.

\begin{corollary} \label{TS26}
Let $A \in {\mathbb C}_{{\rm H}}^{m}$ and  $B \in \mathbb C^{m\times n}$ be given$.$ Then$,$ the following hold$.$
\begin{enumerate}
\item[{\rm(a)}] There exists an $X \in {\mathbb C}_{{\rm H}}^{n}$ such that $BXB^{*} \succcurlyeq  A \succcurlyeq  0$
if and only if ${\mathscr R}(A) \subseteq {\mathscr R}(B),$ namely$,$ the matrix equation
$BY = A$ is consistent$.$  In this case$,$ the general Hermitian solution
can be written as
\begin{equation}
X = B^{\dag}A(B^{\dag})^{*} + UU^{*} + F_BV + V^{*}F_B,
\label{qq231}
\end{equation}
where $U, \, V \in \mathbb C^{n\times n}$ are arbitrary$.$

\item[{\rm(b)}] There exists an $X  \in {\mathbb C}_{{\rm H}}^{n}$ such that
$BXB^{*} \succ  A \succcurlyeq  0$ if and only if  $r(B)= m.$ In this case$,$
the general Hermitian solution can be written as {\rm
(\ref{qq231})}$,$  in which $U, \, V\in \mathbb C^{n\times n}$ are
arbitrary$.$

\item[{\rm(c)}] There exists an $X \in {\mathbb C}_{{\rm H}}^{n}$ such that $BXB^{*}  \preccurlyeq  A \preccurlyeq  0$
if and only if  ${\mathscr R}(A) \subseteq {\mathscr R}(B).$ namely$,$ the matrix equation
$BY = A$ is consistent$.$  In this case$,$ the general Hermitian solution
can be written as
\begin{equation}
X = B^{\dag}A(B^{\dag})^{*}  - UU^{*} + F_AV + V^{*}F_A, \label{qq234}
\end{equation}
where $U,\, V \in \mathbb C^{n\times n}$ are arbitrary$.$

\item[{\rm(d)}] There exists an $X \in \mathbb C^{n\times n}$ such that $BXB^{*}  \prec  A \preccurlyeq  0$
if and only if  $r(B)= m.$ In this case$,$ the general Hermitian solution
 can be written as {\rm (\ref{qq234})}$,$  in which
$U,\ V\in \mathbb C^{n\times n}$ are arbitrary$.$
\end{enumerate}
In particular$,$
\begin{enumerate}
\item[{\rm(e)}] There exists a $0 \preccurlyeq  X \in {\mathbb C}_{{\rm H}}^{n}$ such
that $BXB^{*} \succcurlyeq  A \succcurlyeq  0$ if and only if ${\mathscr R}(A) \subseteq {\mathscr R}(B),$
namely$,$ the matrix equation $BY = A$ is consistent$.$  In this case$,$ the general
positive semi-definite solution  can be written as
\begin{equation}
X = B^{\dag}A(B^{\dag})^{*} + UU^{*}, \label{xx238}
\end{equation}
where $U \in \mathbb C^{n\times n}$ is arbitrary$.$

\item[{\rm(f)}] There exists a $0 \preccurlyeq  X \in {\mathbb C}_{{\rm H}}^{n}$ such
that
\begin{equation}
BXB^{*} \succ  A \succcurlyeq  0
\label{xx239}
\end{equation}
if and only if  $r(B)= m.$ In this case$,$ the general general
positive semi-definite solution of {\rm (\ref{xx239})} can be written
as {\rm (\ref{xx238})}$,$  in which $U \in \mathbb C^{n\times n}$ is
arbitrary$.$
\end{enumerate}
\end{corollary}

\noindent {\bf Proof.} \ Under the condition $A \succcurlyeq  0$,
(\ref{qq27}) is equivalent to $E_BA =0$, i.e., ${\mathscr R}(A)
\subseteq {\mathscr R}(B).$ In this case, (\ref{qq28}) reduces to
(\ref{qq231}).   Also under the condition $A \succcurlyeq  0$,
(\ref{qq210}) is equivalent to $E_B =0$ i.e., $BB^{\dag} = I_m$,
which is further equivalent to $r(B) = m$, as required for (b).
Results (c) and (d) can be shown similarly. Results (e) and (f)
follow from  (a) and (b). \qquad $\Box$

\medskip

Note that the formulas in (\ref{qq28}) and (\ref{qq212b}) are given in closed-form with two independent
parametric matrices. Hence, it is easy to use the two formulas  in the investigation of algebraic
properties of the LMI in (\ref{ww18}) and various problems related to the LMI.
 In Theorem 3.2 below, we shall give the global maximal and minimal
ranks and inertias of the two Hermitian solutions in (\ref{qq28}) and (\ref{qq212b});
in Section 3, we shall use (\ref{qq28}) and (\ref{qq212b}) to solve the
inequality-constrained rank and inertia optimization problems in (\ref{ww19})--(\ref{ww112}).


\section{Ranks and inertias of  $A_1 - B_1XB^{*}_1$ subject to
$B_2XB^{*}_2 \succcurlyeq   A_2$}
\renewcommand{\theequation}{\thesection.\arabic{equation}}
\setcounter{section}{3} \setcounter{equation}{0}

Note that (\ref{qq28}) is in fact a quadratic form, so that (\ref{ww18}) is a quadratic form as well.
 To solve (\ref{ww19})--(\ref{ww112}),  we need the following known results.

\begin{lemma}[\cite{T-mcm,T-laa12}] \label{TS31}
Let  $A \in {\mathbb C}_{{\rm H}}^{m}$ and $B \in \mathbb C^{m \times n}$ be given$,$
and let $M = \left[\!\begin{array}{cccc} A & B \\
B^{*} &0 \end{array}\!\right]\!.$ Then$,$
\begin{align}
\max_{X\in \mathbb C^{n \times k}}\!\!r(\, A + BXX^{*}B^{*} \,) & = \min\{\, r[\, A, \, B\,], \ \ k + r(A) \,\},
\label{243}
\\
\min_{X\in \mathbb C^{n \times k}}\!\!r(\, A + BXX^{*}B^{*} \,) & =  \max\left\{\,  r(A) - k,  \
i_{+}(A) + r[\, A, \, B\,] - i_{+}(M) \right\},
\label{244}
\\
\max_{X\in \mathbb C^{n \times k}}\!\!i_{+}(\, A +BXX^{*}B^{*} \,) & =
\min \left\{\, i_{+}(M), \ \ k +  i_{+}(A) \,\right\},
\label{245}
\\
\min_{X\in \mathbb C^{n \times k}}\!\!i_{+}(\, A +BXX^{*}B^{*} \,)  & =  i_{+}(A),
\label{246}
\\
\max_{X\in \mathbb C^{n \times k}}\!\!i_{-}(\, A +BXX^{*}B^{*} \,) & =  i_{-}(A),
 \label{247}
\\
\min_{X\in \mathbb C^{n \times k}}\!\!i_{-}(\, A +BXX^{*}B^{*} \,) &
\max \{\, i_{-}(A) -k, \ r[\, A, \, B\,] - i_{+}(M) \, \},
 \label{248}
\\
\max_{X\in \mathbb C^{n \times k}}\!\!r(\, A - BXX^{*}B^{*} \,) &   \min\{\, r[\, A, \, B\,], \ \ k + r(A) \,\},
\label{gg326}
\\
\min_{X\in \mathbb C^{n \times k}}\!\!r(\, A - BXX^{*}B^{*} \,)
& =  \max\left\{\,  r(A) - k, \  i_{-}(A) + r[\, A, \, B\,] - i_{-}(M) \right\},
\label{gg327}
\\
\max_{X\in \mathbb C^{n \times k}}\!\!i_{+}(\, A - BXX^{*}B^{*} \,) & =  i_{+}(A),
\label{gg328}
\\
\min_{X\in \mathbb C^{n \times k}}\!\!i_{+}(\, A - BXX^{*}B^{*} \,) & =  \max \{
\, i_{+}(A) - k, \ \  r[\,A, \, B\,] -  i_{-}(M) \, \},
\label{gg329}
\\
\max_{X\in \mathbb C^{n \times k}}\!\!i_{-}(\, A  - BXX^{*}B^{*} \,) & =
\min \left\{\, i_{-}(M), \ \ k +  i_{-}(A) \,\right\},
\label{gg330}
\\
\min_{X\in \mathbb C^{n \times k}}\!\!i_{-}(\, A  - BXX^{*}B^{*} \,) & =  i_{-}(A).
\label{gg331}
\end{align}
In particular$,$
\begin{align}
&  \max_{X\in \mathbb C^{n \times n}}\!\!r(\, A + BXX^{*}B^{*} \,) = r[\, A, \, B\,],
\label{256}
\\
& \min_{X\in \mathbb C^{n \times n}}\!\!r(\, A + BXX^{*}B^{*} \,)   = i_{+}(A) + r[\, A, \, B\,] - i_{+}(M), \ \
\label{256a}
\\
&  \max_{X\in \mathbb C^{n \times n}}\!\!i_{+}(\, A + BXX^{*}B^{*} \,) =
i_{+}(M),
\label{257}
\\
& \min_{X\in \mathbb C^{n \times n}}\!\!i_{+}(\, A + BXX^{*}B^{*} \,)  = i_{+}(A), \ \
\label{258}
\\
& \max_{X\in \mathbb C^{n \times n}}\!\!i_{-}(\, A + BXX^{*}B^{*} \,) = i_{-}(A),
\label{259}
\\
& \min_{X\in \mathbb C^{n \times n}}\!\!i_{-}(\, A + BXX^{*}B^{*} \,)   =r[\, A, \, B\,] - i_{+}(M), \ \
 \label{260}
\\
&  \max_{X\in \mathbb C^{n \times n}}\!\!r(\, A - BXX^{*}B^{*} \,) = r[\, A,
\, B\,],
\label{261}
\\
& \min_{X\in \mathbb C^{n \times n}}r(\, A - BXX^{*}B^{*} \,) = i_{-}(A) + r[\, A, \, B\,] - i_{-}(M),
\label{262}
\\
&  \max_{X\in \mathbb C^{n \times n}}\!\!i_{+}(\, A - BXX^{*}B^{*} \,) =
i_{+}(A),
\label{263}
\\
&  \min_{X\in \mathbb C^{n \times n}}\!\!i_{+}(\, A - BXX^{*}B^{*} \,)  =  r[\,A, \, B\,] -  i_{-}(M),
\label{264}
\\
&  \max_{X\in \mathbb C^{n \times n}}\!\!i_{-}(\, A - BXX^{*}B^{*} \,) =
i_{-}(M),
\label{265}
\\
&  \min_{X\in \mathbb C^{n \times n}}\!\!i_{-}(\, A - BXX^{*}B^{*} \,)  = i_{-}(A).
\label{266}
\end{align}
\end{lemma}

\begin{lemma} [\cite{LT-jota}] \label{TS32}
Let $A_i \in \mathbb C_{{\rm H}}^{m_i}$ and $B_i \in \mathbb
C^{m_i \times n}$ be given, $ i =1,\, 2,$ and assume that the matrix equation
$B_2 XB^{*}_2 = A_2$ has a Hermitian solution$.$ Also let
\begin{eqnarray}
{\cal S} = \{ \, X \in \mathbb C_{{\rm H}}^{n} \ | \  B_2 XB^{*}_2 = A_2 \},  \ \
M =\left[\!\!\begin{array}{ccc} A_1 &  0 &　B_1 \\　0 & -A_2 &
 B_2　\\  B^{*}_1 & B^{*}_2 & 0  \end{array}\!\!\right]\!, \ N = \left[\!\!\begin{array}{ccc} A_1 & B_1 & 0
\\ B_1^{*} & 0 & B_2^{*}
  \end{array}\!\!\right]\!.
 \label{qq31}
\end{eqnarray}
Then$,$
\begin{align}
& \max_{X\in {\cal S}}\!r(\,A_1 - B_1XB_1^{*}\,) = \min\left\{ r[\,A_1,\, B_1\,], \ \
r(M) - 2r(B_2) \right\}\!,
 \label{qq32}
\\
& \min_{X\in {\cal S}}\!r(\,A_1-B_1XB_1^{*}\,) = 2r[\,A_1,\, B_1\,] -
2r(N) + r(M),
 \label{qq33}
\\
& \max_{X\in {\cal S}}\!i_{\pm}(\,A_1-B_1XB_1^{*}\,) =
i_{\pm}(M) - r(B_2),
 \label{qq34}
\\
& \min_{X\in {\cal S}}\!i_{\pm}(\,A_1-B_1XB_1^{*}\,) =
r[\,A_1,\, B_1\,] - r(N)+ i_{\pm}(M).
 \label{qq35}
\end{align}
In consequence$,$ the following hold$.$
\begin{enumerate}
\item[{\rm (a)}] There exists an $X \in \mathbb C_{{\rm H}}^{n}$ such that
$A_1 - B_1XB_1^{*}$ is nonsingular and  $B_2XB_2^{*} = A_2$ if and only if  $r[\,A_1,\, B_1\,] = m_1$
or $r(M) = 2r(B_2) + m_1.$

\item[{\rm (b)}] There exists an $X \in \mathbb C_{{\rm H}}^{n}$ such that $B_1XB_1^{*}=A_1$ and
$B_2XB_2^{*} = A_2$  if and only if
\begin{equation}
{\mathscr R}(A_1)\subseteq {\mathscr R}(B_1),  \ \  {\mathscr R}(A_2)\subseteq {\mathscr R}(B_2), \ \
r(M)= 2r\!\left[\!\!\begin{array}{c} B_1 \\ B_2  \end{array}\!\!\right]\!.
\label{qq36}
\end{equation}

\item[{\rm (c)}] There exists an $X \in \mathbb C_{{\rm H}}^{n}$ such that $B_1XB_1^{*}\prec  A_1$ and
$B_2XB_2^{*} = A_2$ if and only if $i_{+}(M) = r(B_2) + m_1.$

\item[{\rm (d)}] There exists an $X \in \mathbb C_{{\rm H}}^{n}$ such that  $B_1XB_1^{*}\succ  A_1$ and
 $B_2XB_2^{*} = A_2$  if and only if $i_{-}(M) = r(B_2) + m_1.$

\item[{\rm (e)}] There exists an $X \in \mathbb C_{{\rm H}}^{n}$ such that
$B_1XB_1^{*}\preccurlyeq  A_1$ and  $B_2XB_2^{*} = A_2$ if and only if
\begin{equation}
{\mathscr R}(A_2) \subseteq {\mathscr R}(B_2) \ \ and  \ \  i_{+}(M) = r(N)- r[\,A_1,\, B_1\,].
\label{qq37}
\end{equation}

\item[{\rm (f)}] There exists an $X \in \mathbb C_{{\rm H}}^{n}$ such that $B_1XB_1^{*}\succcurlyeq  A_1$
and  $B_2XB_2^{*} = A_2$  if and only if
\begin{equation}
{\mathscr R}(A_2)\subseteq {\mathscr R}(B_2) \ \ and  \ \  i_{-}(M) = r(N) - r[\,A_1,\, B_1\,].
\label{qq38}
\end{equation}
\end{enumerate}
\end{lemma}

We next solve (\ref{ww19})--(\ref{ww112}) for $i = 1$.

\begin{theorem} \label{TS33}
Let $A_i \in \mathbb C_{{\rm H}}^{m_i}$ and $B_i \in \mathbb C^{m_i \times n}$ be given$,$ $ i =1,\, 2,$ and
let $M$ and $N$ be of the forms in {\rm (\ref{qq31})}$.$ Also$,$ assume that there exists an  $X \in {\mathbb
C}_{{\rm H}}^{n}$  such that $B_2XB^{*}_2 \succcurlyeq   A_2,$
 i.e.$,$ $i_{-}\!\left[\!\!\begin{array}{ccc} A_2 & B_2 \\  B_2^{*} & 0 \end{array}\!\!\right] = r[\,A_2,\, B_2\,],
$ and let
\begin{eqnarray}
{\cal T}_1 = \{ \, X \in \mathbb C_{{\rm H}}^{n} \ | \  B_2 XB^{*}_2 \succcurlyeq  A_2 \},  \ \ \ \
 M_1 =\left[\!\!\begin{array}{ccc} A_1 & B_1 \\ B_1^{*} & 0 \end{array}\!\!\right]\!.
 \label{qq39}
\end{eqnarray}
 Then$,$
\begin{align}
& \max_{X\in {\cal T}_1} r(\, A_1 - B_1XB_1^{*}\,)  = r[\, A_1, \, B_1\, ],
\label{qq310}
\\
& \min_{X\in {\cal T}_1} r(\, A_1 - B_1XB_1^{*}\,) =
2r[\,A_1,\, B_1\,] +  i_{-}(M) -  i_{-}(M_1)-  r(N),
\label{qq311}
\\
& \max_{X\in {\cal T}_1} i_{+}(\, A_1 - B_1XB_1^{*}\,) =  i_{+}(M)- r[\, A_2, \, B_2 \,],
\label{qq312}
\\
& \min_{X\in {\cal T}_1} i_{+}(\, A_1 - B_1XB_1^{*}\,) =  r[\,A_1,\, B_1\,]  -
i_{-}(M_1),
\label{qq313}
\\
& \max_{X\in {\cal T}_1} i_{-}(\, A_1 - B_1XB_1^{*}\,) =  i_{-}(M_1),
\label{qq314}
\\
& \min_{X\in {\cal T}_1} i_{-}(\, A_1 - B_1XB_1^{*}\,) =  r[\,A_1,\, B_1 \,] - r(N) + i_{-}(M).
\label{qq315}
\end{align}
In consequence$,$  the following hold$.$
\begin{enumerate}
\item[{\rm (a)}] There exists an $X \in \mathbb C_{{\rm H}}^{n}$ such that
$A_1 - B_1XB_1^{*}$ is nonsingular and $B_2XB_2^{*} \succcurlyeq  A_2$ if and only if  $r[\,A_1,\, B_1\,]
 = m_1.$

\item[{\rm (b)}] The rank of $A_1 - B_1XB_1^{*}$ is nonsingular  for any  $B_2XB_2^{*} \succcurlyeq  A_2$
if and only if $i_{-}(M) =  i_{-}(M_1) + r(N) - m_1.$

\item[{\rm (c)}] There exists an $X \in \mathbb C_{{\rm H}}^{n}$ such that $B_1XB_1^{*}=A_1$ and
$B_2XB_2^{*} \succcurlyeq  A_2$ if and only if ${\mathscr R}(A_1)\subseteq {\mathscr R}(B_1)$  and
$i_{-}(M)= r\!\left[\!\!\begin{array}{c} B_1 \\ B_2  \end{array}\!\!\right]\!.$

\item[{\rm (d)}] There exists an $X \in \mathbb C_{{\rm H}}^{n}$ such that $B_1XB_1^{*}\prec A_1$  and
 $B_2XB_2^{*} \succcurlyeq  A_2$  if and only if $i_{+}(M) = r[\, A_2, \, B_2 \,] +m_1.$

\item[{\rm (e)}] There exists an $X \in \mathbb C_{{\rm H}}^{n}$ such that  $B_1XB_1^{*}\succ A_1$   and
 $B_2XB_2^{*} \succcurlyeq  A_2$ if and only if $i_{-}(M_1) = m_1.$

\item[{\rm (f)}] There exists an $X \in \mathbb C_{{\rm H}}^{n}$ such that  $B_1XB_1^{*} \preccurlyeq  A_1$  and
 $B_2XB_2^{*} \succcurlyeq  A_2$  if and only if $  i_{-}(M)  = r(N) -r[\,A_1,\, B_1 \,].$

\item[{\rm (g)}] There exists an $X \in \mathbb C_{{\rm H}}^{n}$ such that  $B_1XB_1^{*}\succcurlyeq  A_1$  and
 $B_2XB_2^{*} \succcurlyeq  A_2$  if and only if $i_{-}(M_1) = r[\,A_1,\, B_1\,].$

\item[{\rm (h)}] $B_1XB_1^{*} \succcurlyeq  A_1$ holds for all $X \in {\mathbb C}_{{\rm H}}^{n}$ such that
$B_2XB_2^{*} \succcurlyeq  A_2,$ i.e.$,$
$$
\{ \, X \in \mathbb C_{{\rm H}}^{n} \ | \  B_2 XB^{*}_2 \succcurlyeq  A_2 \} \subseteq
\{ \, X \in \mathbb C_{{\rm H}}^{n} \ | \  B_1 XB^{*}_1 \succcurlyeq  A_1 \}
$$
 if and only if  $i_{+}(M)= r[\, A_2, \, B_2 \,].$

\item[{\rm (i)}]  The rank of $A_1 - B_1XB_1^{*}$ is invariant with respect to the Hermitian solution
of $B_2XB_2^{*} \succcurlyeq  A_2$ if and only if $i_{-}(M) = i_{-}(M_1) +  r(N) - r[\,A_1,\, B_1\,].$

\item[{\rm (j)}]  The positive index of the inertia of $A_1 - B_1XB_1^{*}$ is invariant
with respect to the solution of $B_2XB_2^{*} \succcurlyeq  A_2$ if and only if
$i_{+}(M) + i_{-}(M_1)   = r[\,A_1,\, B_1\,] + r[\, A_2, \, B_2 \,].$

\item[{\rm (k)}]  The negative index of the inertia of $A_1 - B_1XB_1^{*}$ is invariant
with respect to the solution of $B_2XB_2^{*} \succcurlyeq  A_2$ if and only if
$i_{-}(M_1) +  r(N)  =  r[\,A_1,\, B_1 \,] + i_{-}(M).$

\item[{\rm (l)}] If  there exist $X_1, \, X_2  \in \mathbb C_{{\rm H}}^{n}$ such that
$B_1X_1B^{*}_1 = A_1$ and  $B_2X_2B^{*}_2 = A_2$ hold$,$ respectively$,$ i.e.$,$
$\R(A_1) \subseteq \R(B_1)$ and $\R(A_2) \subseteq \R(B_2),$ then$,$
\begin{align}
&  \max_{X\in {\cal T}_1} r(\, A_1 - B_1XB_1^{*}\,) = r(B_1),
\\
& \min_{X\in {\cal T}_1} r(\, A_1 - B_1XB_1^{*}\,)  = i_{-}(M) -  r[\, B_1, \, B_2\,],
\label{qq318}
\\
&  \max_{X\in {\cal T}_1} i_{+}(\, A_1 - B_1XB_1^{*}\,)  =  i_{+}(M)- r(B_2),
\\
& \min_{X\in {\cal T}_1} i_{+}(\, A_1 - B_1XB_1^{*}\,) = 0, \ \ \
\label{qq320}
\\
&  \max_{X\in {\cal T}_1} i_{-}(\, A_1 - B_1XB_1^{*}\,)  = r(B_1),
\\
& \min_{X\in {\cal T}_1} i_{-}(\, A_1 - B_1XB_1^{*}\,)  = i_{-}(M) - r[\,B_1,\, B_2 \,].
\label{qq322}
\end{align}
\end{enumerate}
\end{theorem}

\noindent {\bf Proof.} \ From Theorem 2.2(a), the general Hermitian solution of $B_2XB^{*}_2 \succcurlyeq   A_2$
 can be written as
\begin{equation}
X = [\, 0, \, I_n\,]J^{\dag}\!\left[\!\!\begin{array}{ccc} 0 \\
I_n \end{array}\!\!\right] + UU^{*} - F_{B_2}V - V^{*}F_{B_2}, \ \ J =\left[\!\!\begin{array}{ccc} - A_2 & B_2 \\
B^{*}_2 & 0 \end{array}\!\!\right]\!,
\label{qq323}
\end{equation}
where  $U, \, V\in {\mathbb C}^{n \times n}$ are arbitrary$.$
Substituting (\ref{qq323}) into $A_1 - B_1XB_1^{*}$ gives
\begin{eqnarray}
A_1 - B_1XB_1^{*} = \widehat{A} - B_1UU^{*}B_1^{*} + B_1F_{B_2}VB_1^{*} +
B_1V^{*}F_{B_2}B_1^{*},
\label{qq324}
\end{eqnarray}
where $\widehat{A} = A_1 -  [\, 0, \, B_1\,]J^{\dag}\left[\!\!\begin{array}{ccc} 0 \\
B_1^{*} \end{array}\!\!\right].$ In consequence,
\begin{align}
& \max_{X\in {\cal T}_1} r(\, A_1 - B_1XB_1^{*}\,) = \max_{U, \, V\in \mathbb C^{n\times n}}
r(\,\widehat{A}  - B_1UU^{*}B_1^{*} + B_1F_{B_2}VB_1^{*} +
B_1V^{*}F_{B_2}B_1^{*}\,), \label{qq326}
\\
& \min_{X\in {\cal T}_1} r(\, A_1 - B_1XB_1^{*}\,) =  \min_{U, \,
  V\in \mathbb C^{n\times n}}\!\!\!r(\,\widehat{A}  - B_1UU^{*}B_1^{*} + B_1F_{B_2}VB_1^{*} + B_1V^{*}F_{B_2}B_1^{*}\,),
\label{327}
\\
& \max_{X\in {\cal T}_1} i_{\pm}(\, A_1 - B_1XB_1^{*}\,) =  \max_{U, \,
 V\in \mathbb C^{n\times n}} \!\!\!i_{\pm}(\,\widehat{A}  -
  B_1UU^{*}B_1^{*} + B_1F_{B_2}VB_1^{*} + B_1V^{*}F_{B_2}B_1^{*}\,),
\label{328}
\\
& \min_{X\in {\cal T}_1} i_{\pm}(\, A_1 - B_1XB_1^{*}\,) =  \min_{U, \,
 V\in \mathbb C^{n\times n}}\!\!\!i_{\pm}(\,\widehat{A}  - B_1UU^{*}B_1^{*} + B_1F_{B_2}VB_1^{*} + B_1V^{*}F_{B_2}B_1^{*}\,).
\label{329}
\end{align}
Applying (\ref{133})--(\ref{136}) gives
\begin{align}
&  \max_{V\in \mathbb C^{n\times n}}\!\!\!r(\,\widehat{A}  - B_1UU^{*}B_1^{*} +
B_1F_{B_2}VB_1^{*} + B_1V^{*}F_{B_2}B_1^{*}\,) \nb
\\
&   =  \min\left\{ r[\, A_1, \, B_1\, ], \ \ r\!\left[\!\!
\begin{array}{cc} \widehat{A}  - B_1UU^{*}B_1^{*} & B_1F_{B_2}
\\
 F_{B_2}B_1^{*} & 0 \end{array} \!\!\right] \right\},
\label{330}
\\
&  \min_{V\in \mathbb C^{n\times n}}\!\!\!r(\,\widehat{A}  - B_1UU^{*}B_1^{*} +
B_1F_{B_2}VB_1^{*} + B_1V^{*}F_{B_2}B_1^{*}\,)  \nb
\\
&  =  2r[\,A_1,\, B_1\,] + r\!\left[\begin{array}{cc} \widehat{A}  -
B_1UU^{*}B_1^{*} & B_1F_{B_2} \\ F_{B_2}B_1^{*} & 0 \end{array}\right] -
2r\!\left[\begin{array}{cc} A_1  & B_1F_{B_2}
\\ B_1^{*} & 0\end{array}\right],
\label{331}
\\
&  \max_{V\in {\mathbb C}^{n \times n}}\!\!\!i_{\pm}(\,\widehat{A}  -
B_1UU^{*}B_1^{*} + B_1F_{B_2}VB_1^{*} + B_1V^{*}F_{B_2}B_1^{*}\,) =
i_{\pm}\!\left[\!\begin{array}{ccc} \widehat{A}  - B_1UU^{*}B_1^{*}
 & B_1F_{B_2}  \\  F_{B_2}B_1^{*}  & 0    \end{array}\!\right],
 \label{332}
\\
&  \min_{V\in {\mathbb C}^{n \times n}}\!\!\!i_{\pm}(\,\widehat{A}  -
B_1UU^{*}B_1^{*} + B_1F_{B_2}VB_1^{*} + B_1V^{*}F_{B_2}B_1^{*}\,)  \nb
\\
&  =  r[\,A_1,\, B_1\,] + i_{\pm}\!\left[\begin{array}{cc} \widehat{A}
- B_1UU^{*}B_1^{*} & B_1F_{B_2} \\ F_{B_2}B_1^{*} & 0 \end{array}\right] -
r\!\left[\begin{array}{cc} A_1  & B_1F_{B_2}
\\ B_1^{*} & 0\end{array}\right].
\label{333}
\end{align}
Note that
\begin{eqnarray}
\left[\!\!\begin{array}{cc} \widehat{A}  - B_1UU^{*}B_1^{*} & B_1F_{B_2}
\\ F_{B_2}B_1^{*} & 0 \end{array}\!\!\right] =
\left[\!\!\begin{array}{cc} \widehat{A} & B_1F_{B_2} \\ F_{B_2}B_1^{*}
& 0 \end{array}\!\!\right] - \left[\!\!\begin{array}{c} B_1 \\ 0
\end{array}\!\!\right]UU^{*}[\, B_1^{*}, \, 0 \,].
\label{334}
\end{eqnarray}
Then applying (\ref{gg326})--(\ref{gg331}) and simplifying by (\ref{qq111}) and (\ref{116}), we obtain
\begin{align}
& \max_{U\in \mathbb C^{n\times n}} \!\!\!r\!\left(
\left[\!\!\begin{array}{cc} \widehat{A} & B_1F_{B_2} \\ F_{B_2}B_1^{*}
& 0 \end{array}\!\!\right] - \left[\!\!\begin{array}{c} B_1 \\ 0
\end{array}\!\!\right]UU^{*}[\, B_1^{*}, \, 0 \,] \right) \nb
\\
& = r\!\left[\!\!\begin{array}{ccc} \widehat{A} & B_1F_{B_2}  &  B_1
 \\ F_{B_2}B_1^{*} & 0 & 0\end{array}\!\!\right]  = r\!\left[\!\!\begin{array}{cc} A_1 & B_1
 \\ F_{B_2}B_1^{*} & 0 \end{array}\!\!\right] = r\!\left[\!\!\begin{array}{ccc} A_1 & B_1 & 0
 \\ B_1^{*} & 0  & B_2^{*} \end{array}\!\!\right]  - r(B_2) = r(N) - r(B_2),
\label{335}
\\
 &\min_{U\in \mathbb C^{n\times n}} \!\!\!r\!\left(
\left[\!\!\begin{array}{cc} \widehat{A} & B_1F_{B_2} \\ F_{B_2}B_1^{*}
& 0 \end{array}\!\!\right]- \left[\!\!\begin{array}{c} B_1 \\ 0
\end{array}\!\!\right]UU^{*}[\, B_1^{*}, \, 0 \,] \right) \nb
\\
 & =  i_{-}\!\left[\!\!\begin{array}{ccc} \widehat{A} &  B_1F_{B_2} \\
 F_{B_2}B_1^{*} & 0 \end{array}\!\!\right] +
  r\!\left[\!\!\begin{array}{ccc} A_1 & B_1F_{B_2} & B_1\\ F_{B_2}B_1^{*} & 0 & 0  \end{array}\!\!\right]
   - i_{-}\!\left[\!\!\begin{array}{ccc} A_1 & B_1F_{B_2} & B_1 \\ F_{B_2}B_1^{*} & 0 &0
\\
B_1^{*} & 0  &  \end{array}\!\!\right] \nb
\\
 &  =  i_{-}\!\left[\!\!\begin{array}{ccc} \widehat{A} & B_1 & 0 \\
B_1^{*} & 0  & B_2^{*} \\ 0  & B_2 & 0 \end{array}\!\!\right] +  r\!\left[\!\!\begin{array}{ccc} A_1 & B_1 & 0 \\ B_1^{*} & 0
 & B_2^{*} \end{array}\!\!\right]  - 2r(B_2)  - i_{-}\!\left[\!\!\begin{array}{ccc} A_1 & B_1 \\
B_1^{*} & 0 \end{array}\!\!\right],
\label{336}
\\
& \max_{U\in \mathbb C^{n\times n}}\!\!\!i_{+}\!\left(
\left[\!\!\begin{array}{cc} \widehat{A} & B_1F_{B_2}
\\ F_{B_2}B_1^{*} & 0 \end{array}\!\!\right] -
\left[\!\!\begin{array}{c} B_1 \\ 0 \end{array}\!\!\right]UU^{*}[\,
B_1^{*}, \, 0 \,] \right)   =  i_{+}\left[\!\!\begin{array}{cc} \widehat{A} & B_1F_{B_2} \\ F_{B_2}B_1^{*} & 0
\end{array}\!\!\right] \nb
\\
 &  =  i_{+}\left[\!\!\begin{array}{ccc} \widehat{A} & B_1 & 0 \\
B_1^{*} & 0  & B_2^{*} \\ 0  & B_2 & 0 \end{array}\!\!\right]  - r(B_2),
\label{337}
\\
 & \min_{U\in \mathbb C^{n\times n}}\!\!\!i_{+}\!\left(
\left[\!\!\begin{array}{cc} \widehat{A} & B_1F_{B_2} \\ F_{B_2}B_1^{*}
& 0
\end{array}\!\!\right] - \left[\!\!\begin{array}{c} B_1 \\ 0 \end{array}\!\!\right]UU^{*}[\, B_1^{*}, \, 0 \,] \right) \nb
\\
 &  =  r\!\left[\!\!\begin{array}{ccc} \widehat{A} & B_1F_{B_2}  &  B_1
 \\ F_{B_2}B_1^{*} & 0 & 0\end{array}\!\!\right] - i_{-}\!\left[\!\!\begin{array}{ccc} \widehat{A} & B_1F_{B_2} &  B_1
  \\
F_{B_2}B_1^{*} & 0 & 0
\\
B_1^{*} & 0 & 0 \end{array}\!\!\right] \nb
\\
 &  =  r\!\left[\!\!\begin{array}{ccc} A_1 & B_1 & 0 \\
B_1^{*} & 0  & B_2^{*} \end{array}\!\!\right]
 - r(B_2)  - i_{-}\!\left[\!\!\begin{array}{ccc} A_1 & B_1 \\
B_1^{*} & 0 \end{array}\!\!\right],
\label{338}
\\
 & \max_{U\in \mathbb C^{n\times n}}\!\!\!i_{-}\!\left(
\left[\!\!\begin{array}{cc} \widehat{A} & B_1F_{B_2}
\\ F_{B_2}B_1^{*} & 0 \end{array}\!\!\right] -
\left[\!\!\begin{array}{c} B_1 \\ 0 \end{array}\!\!\right]UU^{*}[\,
B_1^{*}, \, 0 \,]\right)   =  i_{-}\!\left[\!\!\begin{array}{ccc} \widehat{A} & B_1F_{B_2} &  B_1 \\
F_{B_2}B_1^{*} & 0 & 0
\\
B_1^{*} & 0 & 0 \end{array}\!\!\right]  \nb
\\
 & = i_{-}\!\left[\!\!\begin{array}{ccc} A_1 & B_1 \\
B_1^{*} & 0 \end{array}\!\!\right],
\label{339}
\\
 &  \min_{U\in \mathbb C^{n\times n}}\!\!\!i_{-}\!\left(
\left[\!\!\begin{array}{cc} \widehat{A} & B_1F_{B_2}
\\ F_{B_2}B_1^{*} & 0 \end{array}\!\!\right] - \left[\!\!\begin{array}{c} B_1 \\ 0
\end{array}\!\!\right]UU^{*}[\, B_1^{*}, \, 0 \,] \right)   =  i_{-}\!\left[\!\!\begin{array}{cc}
\widehat{A} & B_1F_{B_2} \\ F_{B_2}B_1^{*} & 0
\end{array}\!\!\right]  \nb
\\
 & = i_{-}\!\left[\!\!\begin{array}{ccc} \widehat{A} & B_1 & 0 \\
B_1^{*} & 0  & B_2^{*} \\ 0  & B_2 & 0 \end{array}\!\!\right] -r(B_2).
\label{340}
\end{align}
Further, applying  congruence matrix operations gives
\begin{align}
& i_{\pm}\!\left[\!\!\begin{array}{ccc} \widehat{A} & B_1 & 0 \\
B_1^{*} & 0  & B_2^{*} \\ 0  & B_2 & 0 \end{array}\!\!\right]
\\
 & = i_{\pm}\!\left[\!\!\begin{array}{ccc} A_1 -  [\, 0, \, B_1\,]J^{\dag}\left[\!\!\begin{array}{c} 0
 \\
B_1^{*} \end{array}\!\!\right] & B_1 & 0
\\
B_1^{*} & 0  & B_2^{*} \\ 0  & B_2 & 0 \end{array}\!\!\right] = i_{\pm}\!\left[\!\!\begin{array}{ccc} A_1 & B_1 &
\frac{1}{2}[\, 0, \, B_1\,]J^{\dag}\left[\!\!\begin{array}{c} 0 \\ B_2^{*} \end{array}\!\!\right]
\cr
B_1^{*} & 0  & B_2^{*}
\cr
\frac{1}{2} [\, 0, \, B_2\,]J^{\dag}\left[\!\!\begin{array}{c} 0 \\ B_1^{*} \end{array}\!\!\right]
 & B_2 & 0 \end{array}\!\!\right] \nb
\\
 & =  i_{\pm}\!\left[\!\!\begin{array}{ccc} A_1 & B_1 & 0
\cr
B_1^{*} & 0  & B_2^{*}
\cr
0  & B_2 & -[\, 0, \, B_2\,]J^{\dag}\left[\!\!\begin{array}{ccc} 0 \\
B_2^{*} \end{array}\!\!\right] \end{array}\!\!\right]  = i_{\pm}\!\left(\left[\!\!\begin{array}{ccc} A_1 & B_1 & 0
\cr
B_1^{*} & 0  & B_2^{*}
\cr
0  & B_2 & 0 \end{array}\!\!\right] -  \left[\!\!\begin{array}{c} 0
\cr
0
\cr
[\, 0, \, B_2\,]\end{array}\!\!\right]  J^{\dag} \left[\, 0, \ 0, \  \left[\!\!\begin{array}{c} 0 \\
B_2^{*} \end{array}\!\!\right]\, \right] \right) \nb
\end{align}
\begin{align}
 & = i_{\pm}\!\left[\!\!\begin{array}{ccccc} A_1 & B_1 & 0  & 0 & 0
\cr
B_1^{*} & 0  & B_2^{*} & 0 & 0
\cr
0  & B_2 & 0  &  0 &  B_2
\cr
0 & 0  & 0 & - A_2 & B_2
\\
0 & 0 & B_2^{*} & B_2^{*}  & 0 \end{array}\!\!\right] - i_{\pm}(J) \ \ \mbox{(by (\ref{123}))}  = i_{\pm}\!\left[\!\!\begin{array}{ccccc} A_1 & B_1 & 0  & 0 & 0
\cr
B_1^{*} & 0  & 0 &  -B_2^{*} & 0
\cr
0  & 0 & 0  &  0 &  B_2
\cr
0 & -B_2  & 0 & - A_2 & 0
\\
0 & 0 & B_2^{*} & 0  & 0 \end{array}\!\!\right] - i_{\pm}(J)
\nb
\\
 & =  i_{\pm}\!\left[\!\!\begin{array}{ccccc} A_1 & B_1  & 0
\cr
B_1^{*} & 0  &  B_2^{*}
\cr
0 & B_2  &  - A_2
 \end{array}\!\!\right]  + r(B_2) - i_{\pm}(J) \ \ \mbox{(by (\ref{qq19}) and (\ref{qq110}))}   \nb
 \\
&  =   i_{\pm}(M)  + r(B_2) - i_{\mp}\left[\!\!\begin{array}{ccc} A_2 & B_2 \\
B^{*}_2 & 0 \end{array}\!\!\right]  \ \ \mbox{(by (\ref{qq17}) and (\ref{qq18}))}.
 \label{341}
\end{align}
Substituting (\ref{341}) into (\ref{335})--(\ref{340}), and  then (\ref{335})--(\ref{340}) into
(\ref{330})--(\ref{333})   yields (\ref{qq310})--(\ref{qq315}).  Applying Lemma 1.4 to
(\ref{qq310})--(\ref{qq315}) leads to (a)--(g). Result (h) follows from
(\ref{qq310})--(\ref{qq315}). \qquad $\Box$

\medskip

Similarly, we can show the following result.

\begin{theorem} \label{TS34}
Let $A_i \in \mathbb C_{{\rm H}}^{m_i}$ and $B_i \in {\mathbb C}^{m_i \times n}$ be given$,$ $ i =1,\, 2,$
and let $M$ and $N$ be of the forms in {\rm (\ref{qq31})}$.$  Also$,$ assume that there exists an
$X \in \mathbb C_{{\rm H}}^{n}$ such that $B_2XB^{*}_2 \preccurlyeq  A_2,$ i.e.$,$
$i_{+}\!\left[\!\!\begin{array}{ccc} A_2 & B_2 \\  B_2^{*} & 0 \end{array}\!\!\right] = r[\,A_2,\, B_2\,],$
and let
\begin{equation}
M_1 =\left[\!\!\begin{array}{ccc} A_1 & B_1 \\ B_1^{*} & 0 \end{array}\!\!\right]\!,  \ \ \
{\cal T}_2 = \{ \, X \in \mathbb C_{{\rm H}}^{n} \ | \  B_2 XB^{*}_2 \preccurlyeq  A_2 \}.
 \label{342}
\end{equation}
 Then$,$
\begin{align}
& \max_{X\in {\cal T}_2} r(\, A_1 - B_1XB_1^{*}\,) =  r[\, A_1, \, B_1\, ],
\label{343}
\\
& \min_{X\in {\cal T}_2} r(\, A_1 - B_1XB_1^{*}\,)  =
2r[\,A_1,\, B_1\,] +  i_{+}(M) -  i_{+}(M_1) -  r(N),
\label{344}
\\
& \max_{X\in {\cal T}_2} i_{+}(\, A_1 - B_1XB_1^{*}\,) =  i_{+}(M_1),
\label{345}
\\
& \min_{X\in {\cal T}_2} i_{+}(\, A_1 - B_1XB_1^{*}\,) =  r[\,A_1,\, B_1 \,] - r(N) + i_{+}(M),
\label{346}
\\
& \max_{X\in {\cal T}_2} i_{-}(\, A_1 - B_1XB_1^{*}\,) =  i_{-}(M)- r[\, A_2, \, B_2 \,],
\\
& \min_{X\in {\cal T}_2} i_{-}(\, A_1 - B_1XB_1^{*}\,) = r[\,A_1,\, B_1\,]  -
i_{+}(M_1).
\label{347}
\end{align}
In consequence$,$  the following hold$.$
\begin{enumerate}
\item[{\rm (a)}] There exists an $X \in \mathbb C_{{\rm H}}^{n}$ such that
$A_1 - B_1XB_1^{*}$ is nonsingular and $B_2XB_2^{*} \preccurlyeq  A_2$ if and only if  $r[\,A_1,\, B_1\,] = m_1.$

\item[{\rm (b)}] The rank of $A_1 - B_1XB_1^{*}$ is nonsingular  for any  $B_2XB_2^{*} \preccurlyeq  A_2$
if and only if $i_{+}(M) =  i_{+}(M_1) + r(N) - m_1.$

\item[{\rm (c)}] There exists an $X \in \mathbb C_{{\rm H}}^{n}$ such that $B_1XB_1^{*}=A_1$ and
$B_2XB_2^{*} \preccurlyeq  A_2$ if and only if
${\mathscr R}(A_1)\subseteq {\mathscr R}(B_1) \ \ and  \ \
i_{+}(M)= r\!\left[\!\!\begin{array}{c} B_1 \\ B_2  \end{array}\!\!\right]\!.$

\item[{\rm (d)}] There exists an $X \in \mathbb C_{{\rm H}}^{n}$ such that $B_1XB_1^{*}\prec A_1$  and
 $B_2XB_2^{*} \preccurlyeq  A_2$  if and only if  $i_{+}(M_1) = m_1.$

\item[{\rm (e)}] There exists an $X \in \mathbb C_{{\rm H}}^{n}$ such that  $B_1XB_1^{*}\succ A_1$   and
 $B_2XB_2^{*} \preccurlyeq  A_2$ if and only if $i_{-}(M) = r[\, A_2, \, B_2 \,] + m_1.$

\item[{\rm (f)}] There exists an $X \in \mathbb C_{{\rm H}}^{n}$ such that  $B_1XB_1^{*} \preccurlyeq  A_1$  and
 $B_2XB_2^{*} \preccurlyeq  A_2$  if and only if $i_{+}(M_1) = r[\,A_1,\, B_1\,].$

\item[{\rm (g)}] There exists an $X \in \mathbb C_{{\rm H}}^{n}$ such that  $B_1XB_1^{*}\succcurlyeq  A_1$  and
 $B_2XB_2^{*} \preccurlyeq  A_2$  if and only if $i_{+}(M)  = r(N) -r[\,A_1,\, B_1 \,].$

\item[{\rm (h)}] $B_1XB_1^{*} \preccurlyeq  A_1$ holds for all $X \in {\mathbb C}_{{\rm H}}^{n}$ such that
$B_2XB_2^{*} \preccurlyeq  A_2,$ i.e.$,$
$$
\{ \, X \in \mathbb C_{{\rm H}}^{n} \ | \  B_2 XB^{*}_2 \preccurlyeq  A_2 \} \subseteq
\{ \, X \in \mathbb C_{{\rm H}}^{n} \ | \  B_1 XB^{*}_1 \preccurlyeq  A_1 \}
$$
 if and only if  $i_{-}(M)= r[\, A_2, \, B_2 \,].$

\item[{\rm (i)}]  The rank of $A_1 - B_1XB_1^{*}$ is invariant with respect to the Hermitian
solution of $B_2XB_2^{*} \preccurlyeq  A_2$ if and only if
$i_{+}(M) = i_{+}(M_1) +  r(N) - r[\,A_1,\, B_1\,].$

\item[{\rm (j)}]  The positive index of the inertia of $A_1 - B_1XB_1^{*}$ is invariant
with respect to the Hermitian solution of $B_2XB_2^{*} \preccurlyeq  A_2$ if and only if
$i_{+}(M_1) +  r(N)  =  r[\,A_1,\, B_1 \,] + i_{+}(M).$

\item[{\rm (k)}]  The negative index of the inertia of $A_1 - B_1XB_1^{*}$ is invariant
with respect to the Hermitian solution of $B_2XB_2^{*} \preccurlyeq  A_2$ if and only if
$i_{-}(M) + i_{+}(M_1)  = r[\,A_1,\, B_1\,] + r[\, A_2, \, B_2 \,].$
\end{enumerate}
\end{theorem}

Theorem \ref{TS33}(g) shows that the pair of LMIs $B_1XB_1^{*}\succcurlyeq  A_1$ and $B_2XB_2^{*}
\succcurlyeq  A_2$ have a common Hermitian solution if and only if
$$
i_{-}\!\left[\!\!\begin{array}{ccc} A_1 & B_1 \\  B_1^{*} & 0 \end{array}\!\!\right] = r[\,A_1,\, B_1\,], \ \
i_{-}\!\left[\!\!\begin{array}{ccc} A_2 & B_2 \\  B_2^{*} & 0 \end{array}\!\!\right] = r[\,A_2,\, B_2\,],
$$
namely, $B_1XB_1^{*}\succcurlyeq  A_1$ and $B_2XB_2^{*}\succcurlyeq  A_2$ have a Hermitian solution, respectively.
This simple fact makes us to give the following conjecture.

\begin{conjecture} \label{TS44}
The $k$ LMIs
$$
B_1XB_1^{*}\succcurlyeq  A_1,  \ldots, B_kXB_k^{*}\succcurlyeq  A_k
$$
have a common Hermitian solution if and only if each of the $k$ LMIs has a
Hermitian solution.
\end{conjecture}

\section{Global maximal and minimal matrices of $A_1 - B_1XB^{*}_1$ subject to LMIs}
\renewcommand{\theequation}{\thesection.\arabic{equation}}
\setcounter{section}{4}
\setcounter{equation}{0}

In this section, we solve the two LMI-constrained partial ordering optimization problems
 in (\ref{qq16}). Let
\begin{eqnarray}
\phi(X) = A_1 - B_1XB^{*}_1.
\label{357}
\end{eqnarray}
Then,  (\ref{qq16}) is equivalent to finding $X_1,  \, X_2 \in {\mathbb C}^{n}_{{\rm H}}$ such that
\begin{eqnarray}
& & B_2X_1B^{*}_2  \succcurlyeq  A_2  \ \ {\rm and} \ \  \phi(X_1)
\succcurlyeq  \phi(X) \ \ \mbox{for all solutions of} \ \ B_2XB^{*}_2
\succcurlyeq  A_2, \label{358}
\\
& & B_2X_2B^{*}_2  \preccurlyeq  A_2  \ \ {\rm and} \ \  \phi(X_2)
\preccurlyeq  \phi(X) \ \
 \mbox{for all solutions of}  \ \  B_2XB^{*}_2  \preccurlyeq  A_2
\label{359}
\end{eqnarray}
hold, respectively.

\begin{theorem} \label{TS51}
Let $A_i \in \mathbb C_{{\rm H}}^{m_i}$ and $B_i \in {\mathbb C}^{m_i \times n}$ be given for
$i =1,\, 2,$ and assume that $B_2XB^{*}_2 \succcurlyeq  A_2$ is consistent$.$ Then$,$ {\rm (\ref{qq16})}
 has a  solution if and only if
\begin{eqnarray}
\R\left[\!\!\begin{array}{ccc} 0 \\
B^{*}_1 \end{array}\!\!\right] \subseteq \R\left[\!\!\begin{array}{ccc} - A_2 & B_2 \\
B^{*}_2 & 0 \end{array}\!\!\right]\!.
\label{360}
\end{eqnarray}
In this case$,$ the global maximizer  $X_0\in {\mathbb C}_{{\rm H}}^{n}$ of {\rm (\ref{358})} is given by
\begin{eqnarray}
X_0 = [\, 0, \,  I_n\,]\left[\!\!\begin{array}{ccc} - A_2 & B_2 \\
B^{*}_2 & 0 \end{array}\!\!\right]^{\dag}\left[\!\!\begin{array}{ccc} 0 \\
I_n \end{array}\!\!\right],
\label{361}
\end{eqnarray}
and the global maximal matrix in {\rm (\ref{qq16})} can be written as
\begin{eqnarray}
 \max_{\succcurlyeq }\{\,  A_1 - B_1XB^{*}_1  \, | \, B_2XB^{*}_2  \succcurlyeq  A_2, \  X\in {\mathbb C}_{{\rm H}}^{n} \,\} =
 A_1 - [\, 0, \, B_1\,]\left[\!\!\begin{array}{ccc} - A_2 & B_2 \\
B^{*}_2 & 0 \end{array}\!\!\right]^{\dag}\left[\!\!\begin{array}{ccc} 0 \\
B_1^{*} \end{array}\!\!\right]\!,
\label{362}
\end{eqnarray}
which satisfies
\begin{eqnarray}
i_{\pm}\!\left( A_1 - [\, 0, \, B_1\,]\left[\!\!\begin{array}{ccc} - A_2 & B_2 \\
B^{*}_2 & 0 \end{array}\!\!\right]^{\dag}\left[\!\!\begin{array}{ccc} 0 \\
B_1^{*} \end{array}\!\!\right] \right) = i_{\pm}(M) - i_{\pm}\!\left[\!\!\begin{array}{ccc} - A_2 & B_2 \\
B^{*}_2 & 0 \end{array}\!\!\right]\!.
\label{362a}
\end{eqnarray}

\end{theorem}

\noindent {\bf Proof.} \ From Lemma \ref{T11}(i), there exists a Hermitian solution $X_0$ of
$B_2XB^{*}_2  \succcurlyeq  A_2$ such that (\ref{358}) holds if and only if
\begin{eqnarray}
\max_{ B_2XB^{*}_2  \succcurlyeq  A_2} \!\!\!i_{+}[\,\phi(X) - \phi(X_0) \,] =
0. \label{363}
\end{eqnarray}
Note that $\phi(X) - \phi(X_0) = B_1X_0B^{*}_1  - B_1XB^{*}_1.$
Applying (\ref{qq312}) to it gives
\begin{eqnarray}
\max_{ B_2XB^{*}_2  \succcurlyeq  A_2}\!\!\!i_{+}[\,\phi(X) - \phi(X_0) \,] &\!\! = \!\!&
\max_{ B_2XB^{*}_2  \succcurlyeq  A_2} \!\!\!i_{+}(\, B_1X_0B^{*}_1  - B_1XB^{*}_1 \,)  \nb
\\
&\!\! = \!\!&  i_{+}\!\left[\!\!\begin{array}{ccc} B_1X_0B^{*}_1  &  0 &　B_1 \\　0 & -A_2 &
 B_2　\\  B^{*}_1 & B^{*}_2 & 0  \end{array}\!\!\right] - i_{+}\!\left[\!\!\begin{array}{ccc} -A_2 & B_2 \\
B^{*}_2 & 0 \end{array}\!\!\right]\!.
\label{364}
\end{eqnarray}
Substituting (\ref{364}) into (\ref{363}) leads to
$$
i_{+}\!\left[\!\!\begin{array}{ccc} B_1X_0B^{*}_1  &  0 &　B_1 \\　0 & -A_2 &
 B_2　\\  B^{*}_1 & B^{*}_2 & 0  \end{array}\!\!\right] =
 i_{+}\!\left[\!\!\begin{array}{ccc} -A_2 & B_2 \\ B^{*}_2 & 0 \end{array}\!\!\right]\!,
$$
which, by Lemma 1.7(g), is equivalent to
\begin{eqnarray}
\R\left[\!\!\begin{array}{ccc} 0 \\
B^{*}_1 \end{array}\!\!\right] \subseteq \R\left[\!\!\begin{array}{ccc} - A_2 & B_2 \\
B^{*}_2 & 0 \end{array}\!\!\right] \ \ {\rm and} \ \
 B_1X_0B^{*}_1 \preccurlyeq  [\, 0, \,  B_1\,]\left[\!\!\begin{array}{ccc} - A_2 & B_2 \\
B^{*}_2 & 0 \end{array}\!\!\right]^{\dag}\left[\!\!\begin{array}{ccc} 0 \\
B^{*}_1 \end{array}\!\!\right]\!.
\label{364a}
\end{eqnarray}
The general solution of the matrix in (\ref{364a}) can be derived from Corollary \ref{T22a}(a).
Comparing the inequality with (\ref{qq323}), we obtain the special solution in (\ref{361}). \qquad $\Box$

\medskip

The following result can be shown similarly.

\begin{theorem} \label{T45}
Let $A_i \in {\mathbb C}_{{\rm H}}^{m_i}$ and $B_i \in {\mathbb C}^{m_i \times n}$ be given for
$i =1,\, 2,$ and assume that $B_2XB^{*}_2 \preccurlyeq  A_2$ is consistent$.$ Then$,$  the second problem in {\rm (\ref{qq16})}
has a solution if and only if
\begin{eqnarray}
\R\left[\!\!\begin{array}{ccc} 0 \\
B^{*}_1 \end{array}\!\!\right] \subseteq \R\left[\!\!\begin{array}{ccc} - A_2 & B_2 \\
B^{*}_2 & 0 \end{array}\!\!\right]\!.
\label{458}
\end{eqnarray}
In this case$,$ the global minimizer  $X_0\in {\mathbb C}_{{\rm H}}^{n}$ of {\rm (\ref{358})} is given by
\begin{eqnarray}
X_0 = [\, 0, \,  I_n\,]\left[\!\!\begin{array}{ccc} - A_2 & B_2 \\
B^{*}_2 & 0 \end{array}\!\!\right]^{\dag}\left[\!\!\begin{array}{ccc} 0 \\
I_n \end{array}\!\!\right]\!,
\label{459}
\end{eqnarray}
and the global minimal matrix in the second problem of {\rm (\ref{qq16})} can uniquely be written as
\begin{eqnarray}
\min_{\succcurlyeq }\{\,  A_1 - B_1XB^{*}_1  \ | \ B_2XB^{*}_2  \preccurlyeq  A_2 \,\} =
 A_1 - [\, 0, \, B_1\,]\left[\!\!\begin{array}{ccc} - A_2 & B_2 \\
B^{*}_2 & 0 \end{array}\!\!\right]^{\dag}\left[\!\!\begin{array}{ccc} 0 \\
B_1^{*} \end{array}\!\!\right]\!,
\label{460}
\end{eqnarray}
which satisfies
\begin{eqnarray}
i_{\pm}\!\left( A_1 - [\, 0, \, B_1\,]\left[\!\!\begin{array}{ccc} - A_2 & B_2 \\
B^{*}_2 & 0 \end{array}\!\!\right]^{\dag}\left[\!\!\begin{array}{ccc} 0 \\
B_1^{*} \end{array}\!\!\right] \right) = i_{\pm}(M) - i_{\pm}\!\left[\!\!\begin{array}{ccc} - A_2 & B_2 \\
B^{*}_2 & 0 \end{array}\!\!\right]\!.
\label{360a}
\end{eqnarray}
\end{theorem}

It is obvious that the right-hand sides of (\ref{364}) and (\ref{460}) are the Schur complement of
$\left[\!\!\begin{array}{ccc} - A_2 & B_2 \\ B^{*}_2 & 0 \end{array}\!\!\right]$ in the adjoint block
 matrix $M$ in (\ref{qq31}).

\section{Ranks and inertias of the Hermitian solutions of $BXB^{*} \succcurlyeq  A$}
\renewcommand{\theequation}{\thesection.\arabic{equation}}
\setcounter{section}{5}
\setcounter{equation}{0}

Note from (\ref{qq28}) and (\ref{qq212b}) that the Hermitian solutions of (\ref{qq26}) and (\ref{qq211})
are in fact quadratic matrix-valued functions that involve two variable matrices. Hence, we are able to
derive from Lemma 1.11 and Theorem 2.6 a group of formulas for calculating the global maximal and minimal
ranks and inertias of the Hermitian solutions of the two LMIs in (\ref{qq26}) and (\ref{qq211}).

\begin{theorem} \label{TS61}
Let $A \in {\mathbb C}_{{\rm H}}^{m}$ and $B \in \mathbb C^{m \times n}$ be given$.$
\begin{enumerate}
\item[{\rm(a)}]  If there exists an $X \in \mathbb C_{{\rm H}}^{n}$  such that $BXB^{*} \succcurlyeq  A$  holds$,$
then$,$
\begin{align}
& \max_{BXB^{*} \succcurlyeq  A} r(X)  = \max_{BXB^{*} \succcurlyeq  A} i_{+}(X)  = n,
\label{ss334}
\\
&  \min_{BXB^{*} \succcurlyeq  A} r(X)
 = \min_{BXB^{*} \succcurlyeq  A} i_{+}(X)  =i_{+}(A),
\label{ss335}
\\
& \max_{BXB^{*} \succcurlyeq  A} i_{-}(X)   =  n + i_{-}(A) - r[\, A, \,  B \,],
\label{ss336}
\\
& \min_{BXB^{*} \succcurlyeq  A} i_{-}(X)  = 0.
\label{ss336a}
\end{align}
In consequence$,$ the following hold$.$
\begin{enumerate}
\item[{\rm(i)}]  The matrix $\widehat{X}$ in {\rm (\ref{qq226})} is a solution that
satisfies the second equality in {\rm (\ref{ss334})}$.$

\item[{\rm(ii)}] $BXB^{*} \succcurlyeq  A$ always has a solution $X \succ 0.$

\item[{\rm(iii)}] All solutions of $BXB^{*} \succcurlyeq  A$ satisfy $X \succ 0$ if and only if $i_{+}(A) =n.$

\item[{\rm(iv)}] $X  = 0$ is a solution of $BXB^{*} \succcurlyeq  A$ $\Leftrightarrow$
 $BXB^{*} \succcurlyeq  A$ has a solution $X \preccurlyeq  0$ $\Leftrightarrow$  $A \preccurlyeq  0.$

\item[{\rm(v)}] $BXB^{*} \succcurlyeq  A$ has a solution $X \prec 0$ if and only if $A\prec 0.$

\item[{\rm(vi)}] All solutions of $BXB^{*} \succcurlyeq  A$ satisfy $X \succcurlyeq  0$ if and only if
 $r(B) = n$ and $ r[\, A, \,  B \,] = i_{-}(A) + n.$
\end{enumerate}

\item[{\rm(b)}]  If there exists an $X \in \mathbb C_{{\rm H}}^{n}$  such that $BXB^{*} \preccurlyeq  A,$  then$,$
\begin{align}
& \max_{BXB^{*} \preccurlyeq  A} r(X)  = \max_{BXB^{*} \preccurlyeq  A} i_{-}(X)  = n,
\label{ss338}
\\
&  \min_{BXB^{*} \preccurlyeq  A} r(X)  = \min_{BXB^{*} \preccurlyeq  A} i_{-}(X)  =i_{-}(A),
\label{ss338a}
\\
&  \max_{BXB^{*} \preccurlyeq  A} i_{+}(X)   =  n + i_{+}(A) - r[\, A, \,  B \,],
\\
& \min_{BXB^{*} \preccurlyeq  A} i_{+}(X)   = 0.
\label{ss340}
\end{align}
\end{enumerate}
In consequence$,$  the following hold$.$
\begin{enumerate}
\item[{\rm(i)}]  The matrix $\widehat{X}$ in {\rm (\ref{qq226})} is a solution that
satisfies the second equality in {\rm (\ref{ss338})}$.$

\item[{\rm(ii)}] $BXB^{*} \preccurlyeq  A$ always has a solution $X \prec 0.$

\item[{\rm(iii)}] All solutions of $BXB^{*} \preccurlyeq  A$ satisfy $X \prec 0$ if and only if $i_{-}(A) =n.$

\item[{\rm(iv)}] $X  = 0$ is a solution of $BXB^{*} \preccurlyeq  A$ $\Leftrightarrow$
 $BXB^{*} \preccurlyeq  A$ has a solution $X \succcurlyeq  0$ $\Leftrightarrow$  $A \succcurlyeq  0.$

\item[{\rm(v)}] $BXB^{*} \preccurlyeq  A$ has a solution $X \succ 0$ if and only if $A\succ 0.$

\item[{\rm(vi)}] All solutions of $BXB^{*} \preccurlyeq  A$ satisfy $X \preccurlyeq  0$ if and only if
 $r(B) = n$ and $ r[\, A, \,  B \,] = i_{+}(A) + n.$
\end{enumerate}

\end{theorem}

Theorem \ref{TS61} shows the ranks and inertias of Hermitian solutions of the two simple
LMIs $BXB^{*} \succcurlyeq  A$ and $BXB^{*} \preccurlyeq  A$ may have different values. In general,
solutions of LMEs and LMIs with low ranks or inertias are objects of particular interest
in the investigations of LMEs and LMIs and their applications; some recent work on
this topic can be found, e.g., in \cite{AHZ,RFP}.

We now turn our attention to the ranks and inertias of submatrices in Hermitian solutions of
the matrix equation $BXB^{*} = A$ and the inequality $BXB^{*} \succcurlyeq  A$ $(BXB^{*} \preccurlyeq  A)$.
 Rewrite $BXB^{*} \succcurlyeq   A$ as
\begin{equation}
  [\, B_1,  \,  B_2 \,]\! \left[\!\! \begin{array}{cc}  X_1  & X_2
\\ X_2^{*} & X_3 \end{array} \!\!\right]\!
\left[\!\! \begin{array}{c} B^{*}_1 \\ B^{*}_2 \end{array} \!\!\right] =  A,
\label{51}
\end{equation}
where $B_1 \in {\mathbb C}^{m \times n_1},$ $B_2 \in \mathbb
C^{m \times n_2}$, $ X_1 \in {\mathbb C}_{{\rm H}}^{n_1}, \ X_2 \in
\mathbb C^{n_1 \times n_2}$ and $X_3 \in  \mathbb C_{{\rm H}}^{n_2}$
 with $n_1 + n_2 = n$.  Note that the submatrices $X_1,  \, X_2,  \, X_3$ in (\ref{51})
 can be rewritten as
\begin{eqnarray}
 X_1 = P_1XP^{*}_1, \ \ X_2 = P_1XP^{*}_2, \ \ X_3= P_2XP^{*}_2,
\label{52}
\end{eqnarray}
where $ P_1 = [\, I_{n_1},  \,  0 \,]$ and $P_2 =[ \, 0, \, I_{n_2} \, ].$
For convenience, we adopt the following notation for the collections of the
submatrices $X_1$ and $X_3$  in (\ref{51}):
\begin{align}
&  {\cal S}_1  =  \left\{ X_1 \in \mathbb C_{{\rm H}}^{n_1} \ \left| \ [\, B_1,\, B_2 \,]
\!\left[\!\! \begin{array}{cc}  X_1  & X_2 \\ X^{*}_2 & X_3
\end{array} \!\!\right]\!\left[\!\! \begin{array}{c} B^{*}_1 \\
B^{*}_2 \end{array} \!\!\right] = A \right.  \right\} =
\left\{ X_1 = P_1XP^{*}_1 \ | \ BXB^{*} =  A, \  X \in \mathbb C_{{\rm H}}^{n} \right\},
\label{53}
\\
&  {\cal S}_3  =  \left\{ X_3 \in \mathbb C_{{\rm H}}^{n_2} \ \left| \ [\, B_1,\, B_2 \,]
\!\left[\!\! \begin{array}{cc}  X_1  & X_2 \\ X^{*}_2 & X_3
\end{array} \!\!\right]\!\left[\!\! \begin{array}{c} B^{*}_1 \\
B^{*}_2 \end{array} \!\!\right] = A \right.  \right\} = \left\{ X_3 = P_2XP^{*}_2 \ | \ BXB^{*} =  A, \
X \in \mathbb C_{{\rm H}}^{n} \right\}.
\label{54}
\end{align}
Applying Lemma \ref{TS32} to (\ref{53}) and (\ref{54}) and simplifying, we obtain the global maximal and minimal ranks
and inertias of the submatrices $X_1$ and $X_3$ in a Hermitian solution of (\ref{51}) and their consequences
as follows. The details of the proof are omitted.

\begin{theorem} \label{TS62}
Assume that the matrix equation {\rm (\ref{51})} is consistent$,$ and let ${\cal S}_1$ and
${\cal S}_3$ be of the forms in {\rm (\ref{53})} and {\rm (\ref{54})}$.$ Then$,$
\begin{align}
&  \max_{X_1 \in {\cal S}_1}r(X_1)  = \min\!\left\{ n_1,  \
r\!\left[\!\! \begin{array}{cc}  A  & B_2 \\
  B^{*}_2 & 0  \end{array} \!\!\right] - 2r(B) + 2n_1 \right\}\!,
\\
&   \min_{X_1 \in {\cal S}_1}r(X_1)  = r\!\left[\!\! \begin{array}{cc}  A & B_2 \\
B^{*}_2 & 0  \end{array} \!\!\right] - 2r(B_2),
\label{56}
\\
& \max_{X_1 \in {\cal S}_1}i_{\pm}(X_1)   =i_{\pm}\!\left[\!\! \begin{array}{cc}  C  & B_2 \\
  A^{*}_2 & 0  \end{array} \!\!\right] - r(A) + n_1,
  \\
& \min_{X_1 \in {\cal S}_1}i_{\pm}(X_1)  = i_{\pm}\!\left[\!\! \begin{array}{cc}  C & B_2 \\
A^{*}_2 & 0  \end{array} \!\!\right] - r(B_2),
\label{58}
\\
&  \max_{X_3 \in {\cal S}_3}r(X_3)  = \min \!\left\{ n_2, \
r\!\left[\!\! \begin{array}{cc}  A  & B_1 \\
  B^{*}_1 & 0  \end{array} \!\!\right] - 2r(B) + 2n_2 \right\}\!,
\\
&  \min_{X_3 \in {\cal S}_3}r(X_3) = r\!\left[\!\! \begin{array}{cc}  A & B_1 \\
B^{*}_1 & 0  \end{array} \!\!\right] - 2r(B_1),
\label{510}
\\
 & \max_{X_3 \in {\cal S}_3}i_{\pm}(X_3) =i_{\pm}\!\left[\!\! \begin{array}{cc}  C  & B_1 \\
  A^{*}_1 & 0  \end{array} \!\!\right] - r(A) + n_2,
  \\
& \min_{X_3 \in {\cal S}_3}i_{\pm}(X_3)  = i_{\pm}\!\left[\!\! \begin{array}{cc}  C & B_1 \\
A^{*}_1 & 0  \end{array} \!\!\right] - r(B_1).
\label{512}
\end{align}
In consequence$,$  the following hold$.$
\begin{enumerate}
\item[{\rm (c)}]  Eq.~{\rm (\ref{51})} has a Hermitian solution in which $X_1$ is nonsingular
if and only if
$
r\!\left[\!\! \begin{array}{cc}  C & B_2 \\ A^{*}_2 & 0  \end{array} \!\!\right] \geqslant   2r(A) - n_1.
$

\item[{\rm (d)}] The submatrix $X_1$ in any Hermitian solution of {\rm (\ref{51})} is nonsingular
if and only if $r\!\left[\!\! \begin{array}{cc}  C & B_2 \\ A^{*}_2 & 0  \end{array} \!\!\right]=  2r(B_2) + n_1.$

\item[{\rm (e)}]  Eq.~{\rm (\ref{51})} has a Hermitian solution in which  $X_1 =0$
if and only if $r\!\left[\!\! \begin{array}{cc}  C & B_2 \\ A^{*}_2 & 0  \end{array} \!\!\right] = 2r(B_2).$

\item[{\rm (f)}]  The submatrix $X_1$ in any Hermitian solution of {\rm (\ref{51})} satisfies $X_1 =0$
if and only if $r\!\left[\!\! \begin{array}{cc}  C & B_2 \\ A^{*}_2 & 0  \end{array} \!\!\right]
=  2r(A) - 2n_1.$

\item[{\rm (g)}]  Eq.~{\rm (\ref{52})} has a Hermitian solution in which $ X_1 \succ  0$ $(X_1 \prec 0)$
if and only if
$$
i_{+}\!\left[\!\! \begin{array}{cc}  C & B_2 \\ A^{*}_2 & 0  \end{array} \!\!\right] =  r(A)
 \ \left(i_{-}\!\left[\!\! \begin{array}{cc}  C  & B_2 \\ A^{*}_2 & 0  \end{array} \!\!\right] =  r(A)
 \right)\!.
$$

\item[{\rm (h)}]   The submatrix $X_1$ in any Hermitian solution of {\rm (\ref{51})} satisfies  $ X_1 \succ  0$ $(X_1 \prec 0)$
 if and only if
$$
i_{+}\!\left[\!\! \begin{array}{cc}  C  & B_2 \\ A^{*}_2 & 0  \end{array} \!\!\right] = n_1 + r(B_2)
\ \left(
i_{-}\!\left[\!\! \begin{array}{cc}  C  & B_2 \\ A^{*}_2 & 0  \end{array} \!\!\right] =
 n_1 + r(B_2) \right)\!.
$$

\item[{\rm (i)}] Eq.~{\rm (\ref{51})} has a Hermitian solution satisfying $ X_1 \succcurlyeq  0$ $(X_1 \preccurlyeq  0)$
if and only if
$$
i_{-}\!\left[\!\! \begin{array}{cc}  C  & B_2 \\ A^{*}_2 & 0
\end{array} \!\!\right] =  r(B_2) \ \ \left(i_{+}\!\left[\!\!
\begin{array}{cc}  C  & B_2 \\ A^{*}_2 & 0  \end{array} \!\!\right]
=  r(B_2) \right)\!.
$$

\item[{\rm (j)}]  The submatrix $X_1$ in any Hermitian solution of {\rm (\ref{51})} satisfies
$X_1 \succcurlyeq  0$ $(X_1 \preccurlyeq  0)$ if and only if
$$
i_{-}\!\left[\!\! \begin{array}{cc}  C & B_2 \\ A^{*}_2 & 0  \end{array} \!\!\right] =  r(A) -n_1
 \ \left( i_{+}\!\left[\!\! \begin{array}{cc}  C  & B_2 \\ A^{*}_2 & 0  \end{array} \!\!\right]
 = r(A) - n_1 \right)\!.
$$

\item[{\rm (k)}] The positive signature of $X_1$ in {\rm (\ref{51})} is invariant $\Leftrightarrow$
the negative signature of $X_1$ {\rm (\ref{51})} is invariant
$\Leftrightarrow$ ${\mathscr R}(B_1) \cap {\mathscr R}(B_2) = \{ 0\}$ and $r(B_1) = n_1.$
\end{enumerate}
\end{theorem}

Replacing the equality sign in (\ref{51}) with inequality signs gives
\begin{align}
 &   [\, B_1,  \,  B_2 \,]\! \left[\!\! \begin{array}{cc}  X_1  & X_2
\\ X_2^{*} & X_3 \end{array} \!\!\right]\!
\left[\!\! \begin{array}{c} B^{*}_1 \\ B^{*}_2 \end{array} \!\!\right] \succcurlyeq  A,
\label{513}
\\
& [\, B_1,  \,  B_2 \,]\! \left[\!\! \begin{array}{cc}  X_1  & X_2
\\ X_2^{*} & X_3 \end{array} \!\!\right]\!
\left[\!\! \begin{array}{c} B^{*}_1 \\ B^{*}_2 \end{array} \!\!\right] \preccurlyeq   A.
\label{514}
\end{align}
Also let
 \begin{align}
&  {\cal U}_1  =  \left\{ X_1 \in \mathbb C_{{\rm H}}^{n_1} \ \left| \ [\, B_1,\, B_2 \,]
\!\left[\!\! \begin{array}{cc}  X_1  & X_2 \\ X^{*}_2 & X_3
\end{array} \!\!\right]\!\left[\!\! \begin{array}{c} B^{*}_1 \\
B^{*}_2 \end{array} \!\!\right] \succcurlyeq  A  \right.  \right\} = \left\{ X_1 = P_1XP^{*}_1 \ | \ BXB^{*} \succcurlyeq   A,
X \in \mathbb C_{{\rm H}}^{n}
\right\},
\label{515}
\\
&  {\cal U}_3 =  \left\{ X_3 \in \mathbb C_{{\rm H}}^{n_2} \ \left| \ [\, B_1,\, B_2 \,]
\!\left[\!\! \begin{array}{cc}  X_1  & X_2 \\ X^{*}_2 & X_3
\end{array} \!\!\right]\!\left[\!\! \begin{array}{c} B^{*}_1 \\
B^{*}_2 \end{array} \!\!\right] \succcurlyeq   A  \right.  \right\} =
\left\{ X_2 = P_2XP^{*}_2 \ | \ BXB^{*} \succcurlyeq   A,  \ X \in \mathbb C_{{\rm H}}^{n} \right\},
\label{516}
\\
& {\cal V}_1 = \left\{ X_1 \in \mathbb C_{{\rm H}}^{n_1} \ \left| \ [\, B_1,\, B_2 \,]
\!\left[\!\! \begin{array}{cc}  X_1  & X_2 \\ X^{*}_2 & X_3
\end{array} \!\!\right]\!\left[\!\! \begin{array}{c} B^{*}_1 \\
B^{*}_2 \end{array} \!\!\right] \preccurlyeq   A \right.  \right\} =
\left\{ X_1 =P_1XP^{*}_1 \ | \ BXB^{*} \preccurlyeq   A, \ X \in \mathbb C_{{\rm H}}^{n} \right\},
\label{517}
\\
&  {\cal V}_3 = \left\{ X_3 \in \mathbb C_{{\rm H}}^{n_2} \ \left| \ [\, B_1,\, B_2 \,]
\!\left[\!\! \begin{array}{cc}  X_1  & X_2 \\ X^{*}_2 & X_3
\end{array} \!\!\right]\!\left[\!\! \begin{array}{c} B^{*}_1 \\
B^{*}_2 \end{array} \!\!\right] \preccurlyeq   A \right.  \right\} =
\left\{ X_3 = P_2XP^{*}_2 \ | \ BXB^{*} \preccurlyeq   A,  \ X \in \mathbb C_{{\rm H}}^{n}\right\}.
\label{518}
\end{align}
Applying  Theorems 3.3 and 3.4 to (\ref{513})--(\ref{518}) and simplifying, we obtain
the global maximal and minimal ranks and inertias of the submatrices $X_1$ and $X_3$ in Hermitian solutions
of (\ref{513}) and (\ref{514})  as follows. The details of the proof are omitted.

\begin{theorem} \label{T53}
Assume that the matrix inequality in {\rm (\ref{513})} is consistent$,$ and let ${\cal U}_1$ and
${\cal U}_3$ be of the forms in {\rm (\ref{515})} and {\rm (\ref{516})}$.$ Then$,$
\begin{eqnarray}
 \max_{X_1 \in {\cal U}_1}r(X_1) &\!\! = \!\!&  \max_{X_1 \in {\cal U}_1}i_{+}(X_1)  = n_1,
\label{519}
\\
\min_{X_1 \in {\cal U}_1}r(X_1) &\!\! = \!\!&  n_1 +  i_{+}\!\left[\!\! \begin{array}{cc}  A & B_2 \\
B^{*}_2 & 0  \end{array} \!\!\right] - r(B_1) - r(B_2),
\label{520}
\\
\min_{X_1 \in {\cal U}_1}i_{+}(X_1) &\!\! = \!\!&  i_{+}\!\left[\!\! \begin{array}{cc}  A & B_2 \\
B^{*}_2 & 0  \end{array} \!\!\right] - r(B_2),
\label{522}
\\
\max_{X_1 \in {\cal U}_1}i_{-}(X_1) &\!\! = \!\!&  n_1 +  i_{-}\!\left[\!\! \begin{array}{cc}  A & B_2 \\
B^{*}_2 & 0  \end{array} \!\!\right] - r[\, A, \, B\,],
\label{523}
\\
\min_{X_1 \in {\cal U}_1}i_{-}(X_1) &\!\! = \!\!&  0,
\label{524}
\\
 \max_{X_1 \in {\cal U}_3}r(X_3) &\!\! = \!\!&  \max_{X_3 \in {\cal U}_3}i_{+}(X_3)  = n_2,
\label{525}
\\
\min_{X_1 \in {\cal U}_3}r(X_3) &\!\! = \!\!&  n_2 +  i_{+}\!\left[\!\! \begin{array}{cc}  A & B_1 \\
B^{*}_1 & 0  \end{array} \!\!\right] - r(B_1) - r(B_2),
\label{526}
\\
\min_{X_3 \in {\cal U}_3}i_{+}(X_3) &\!\! = \!\!&  i_{+}\!\left[\!\! \begin{array}{cc}  A & B_1 \\
B^{*}_1 & 0  \end{array} \!\!\right] - r(B_1),
\label{528}
\\
\max_{X_3 \in {\cal U}_3}i_{-}(X_3) &\!\! = \!\!&  n_2 +  i_{-}\!\left[\!\! \begin{array}{cc}  A & B_1 \\
B^{*}_1 & 0  \end{array} \!\!\right] - r[\, A, \, B\,],
\label{529}
\\
\min_{X_3 \in {\cal U}_3}i_{-}(X_3) &\!\! = \!\!&  0.
\label{530}
\end{eqnarray}
\end{theorem}

\begin{theorem} \label{T54}
Assume that the matrix inequality in {\rm (\ref{514})} is consistent$,$ and let ${\cal V}_1$ and
${\cal V}_3$ be of the forms in {\rm (\ref{517})} and {\rm (\ref{518})}$.$ Then$,$
\begin{align}
 \max_{X_1 \in {\cal V}_1}r(X_1) &  =  \max_{X_1 \in {\cal V}_1}i_{-}(X_1)  = n_1,
\label{531}
\\
\min_{X_1 \in {\cal V}_1}r(X_1) & =   n_1 +  i_{-}\!\left[\!\! \begin{array}{cc}  A & B_2 \\
B^{*}_2 & 0  \end{array} \!\!\right] - r(B_1) - r(B_2),
\label{532}
\\
\max_{X_1 \in {\cal V}_1}i_{+}(X_1) & =   n_1 +  i_{+}\!\left[\!\! \begin{array}{cc}  A & B_2 \\
B^{*}_2 & 0  \end{array} \!\!\right] - r[\, A, \, B\,],
\label{533}
\\
\min_{X_1 \in {\cal V}_1}i_{+}(X_1) & =   0,
\label{534}
\\
\min_{X_1 \in {\cal V}_1}i_{-}(X_1) & =  i_{-}\!\left[\!\! \begin{array}{cc}  A & B_2 \\
B^{*}_2 & 0  \end{array} \!\!\right] - r(B_2),
\label{536}
\\
 \max_{X_3 \in {\cal V}_3}r(X_3) & =   \max_{X_3 \in {\cal V}_3}i_{-}(X_3)  =  n_2,
\label{537}
\\
\min_{X_3 \in {\cal V}_3}r(X_3) & = n_2 +  i_{-}\!\left[\!\! \begin{array}{cc}  A & B_1 \\
B^{*}_1 & 0  \end{array} \!\!\right] - r(B_1) - r(B_2),
\label{538}
\\
\max_{X_3 \in {\cal V}_3}i_{+}(X_3) &=  n_2 +  i_{+}\!\left[\!\! \begin{array}{cc}  A & B_1 \\
B^{*}_1 & 0  \end{array} \!\!\right] - r[\, A, \, B\,],
\label{539}
\\
\min_{X_3 \in {\cal V}_3}i_{+}(X_3) & =  0,
\label{540}
\\
\min_{X_3 \in {\cal V}_3}i_{-}(X_3) & = i_{-}\!\left[\!\! \begin{array}{cc}  A & B_1 \\
B^{*}_1 & 0  \end{array} \!\!\right] - r(B_1).
\label{542}
\end{align}
\end{theorem}

A further work is to give the extremal ranks and inertias of the $A_1 - B_1XB^{*}_1$
subject to the common Hermitian solution of the $k-1$ consistent LMIs
$$
B_2XB^{*}_2 \succcurlyeq  A_2, \ldots, B_kXB_k^{*} \succcurlyeq  A_k,
$$
and to establish necessary and sufficient condition for the set of LMIs
$$
B_1XB^{*}_1 \succcurlyeq  A_1, \  B_2XB^{*}_2 \succcurlyeq  A_2, \ldots, B_kXB_k^{*} \succcurlyeq  A_k,
$$
to have a common Hermitian solution.

Finally, it should be pointed out that the rank and inertia of a matrix, as two simplest concepts in linear
algebra, are also one of the richest fields in mathematics that admit ten thousands of analytical formulas.\\

\noindent {\bf Acknowledgements} \  This work was supported by National Natural Science Foundation of
China (Grant No. 11271384).

\end{document}